\documentclass[10pt]{article}
\usepackage{amsmath,amssymb}
\newtheorem{theorem}{Theorem}[section]
\newtheorem{lemma}{Lemma}[section]
\newtheorem{remark}{Remark}[section]

\newtheorem{proposition}{Proposition}[section]

\newcommand{\proof}{\noindent \textbf {Proof.} }

\newcommand{\qed}{\hfill \hfill $\Box$}
\newcommand\step[1]{\noindent\emph{#1}}

\date{To appear in \emph{Nonlinearity}}
\begin{document}
\title {Stability of traveling wavefronts in discrete   reaction-diffusion equations with nonlocal delay effects   \thanks{This work was partially supported by the National Natural     Science Foundation of P.R.\ China (Grant No. 11271115) and and the UK's EPSRC (EP/K027743/1).}}
\author {Shangjiang Guo\thanks{E-mail: shangjguo@hnu.edu.cn}\\
  College of Mathematics and Econometrics, Hunan
  University\\
  Changsha, Hunan 410082, People's Republic of China\\
  Johannes Zimmer
  \\
  Department of Mathematical Sciences, University of Bath, \\
  Bath BA2 7AY, United Kingdom }
 \maketitle
\begin{abstract}

  This paper deals with traveling wavefronts for temporally delayed,   spatially discrete reaction-diffusion equations. Using a combination   of the weighted energy method and the Green function technique, we   prove that all noncritical wavefronts are globally exponentially   stable, and critical wavefronts are globally algebraically stable   when the initial perturbations around the wavefront decay to zero   exponentially near minus infinity regardless of the magnitude of   time delay.

  {\bf Keywords.}  Travelling waves, time delay, global stability,   Fisher-KPP equation, weighted energy, Green functions

  {\bf AMS subject classifications.} 35K57, 34K20, 92D25.
\end{abstract}

\section{Introduction}
\label{sec:Introduction}
Traveling wavefront solutions play an important role in the description of the long-term behaviour of solutions to initial value problems in reaction-diffusion equations, both in the spatially continuous case and in spatially discrete situations. Such solutions are also of interest in their own right, for example to understand transitions between different states of a physical system, propagation of patterns, and domain invasion of species in population biology (see, e.g.,~\cite{Bates1,Bates2,Cahn,Chow,So}). 
Here we present a stability analysis for traveling wavefronts of the equation
\begin{equation}
  \label{eq}
  u_t(x,t)= d\cdot\Delta_1u(x,t)+f
  \left(u(x,t),(h*u)(x,t-\tau)\right),
\end{equation}
where $x\in\mathbb{R}$, $d>0$, $\tau\geq 0$, $\Delta_1u(x,t)=u(x+1,t)-2u(x,t)+u(x-1,t)$ and
\begin{equation*}
  (h*u)(x,t)=\int_{\mathbb{R}}h(x-y)u(y,t)dy.
\end{equation*}

One interpretation is that equation (\ref{eq}) is a spatially discrete version of a time-delayed reaction-diffusion equations with nonlocal nonlinearity,
\begin{equation}\label{eqcon1}
  u_t(x,t)= du_{xx}(x,t)+f
  \left(u(x,t),(h*u)(x,t-\tau)\right).
\end{equation}
Such reaction-diffusion equations with nonlocal terms have been proposed by Britton~\cite{Britton1,Britton2} in a biological context for population dynamics.  In recent years, spatially non-local differential equations such as (\ref{eqcon1}) have attracted significant attention (see, e.g., \cite{Gourley,So,Thieme, Wang, Ou}). One can expect traveling wave solutions to play an important role in understanding the dynamics of~\eqref{eq} and~\eqref{eqcon1}, as for many other equations. However, one soon sees that the structure of such solutions for system~\eqref{eq} is richer and much more complex than for~\eqref{eqcon1}, as it is a nonlocal equation (see the discussion following~\eqref{eq3}). The study of traveling wave solutions of~\eqref{eqcon1} leads to the second order ordinary differential equation (ODE)
\begin{equation}
  \label{eq4}
      -c\phi'(\xi)+d\phi''(\xi)+f\left(\phi(\xi),
      \int_{\mathbb{R}}h(y)S(\varphi(\xi-y-c\tau))dy\right)=0,\quad \xi\in \mathbb{R},\\
\end{equation}
for $\varphi \colon \mathbb{R}\to \mathbb{R}$, with appropriate boundary conditions for $\varphi(\pm\infty)$ (for convenience, we write $\phi(-\infty)$ and $\phi(\infty)$ as abbreviations for $\lim_{\xi\to-\infty}\phi(\xi)$ and $\lim_{\xi\to \infty}\phi(\xi)$, respectively).  Under some monostable assumptions, Wang \textit{et   al}. \cite{Wang} investigated the existence, uniqueness, and global asymptotical stability of traveling wave fronts for~\eqref{eq4}. We also refer to So \emph{et al}.~\cite{So} for more details and some specific forms of $f$, obtained from integrating along characteristics of a structured population model, an idea from the work of Smith and Thieme~\cite{Smith}. See also~\cite{So} for a similar model and~\cite{Gourley0} for a survey on the history and the current status of the study of reaction diffusion equations with non-local delayed interactions.  In particular, when $f(u,v)= v(1-u)$ and $h(u)=\delta(u)$, the equation (\ref{eqcon1}) is delayed \emph{Fisher's equation}~\cite{Fisher} or \emph{KPP   equation}~\cite{Kolmogorov}, which arises in the study of gene development or population dynamics. When $f(u,v)=-au+b(v)$ and $h(u)=\delta(u)$, the equation (\ref{eqcon1}) is the local \emph{Nicholson's blowflies equation} and has been investigated in~\cite{Gourley,Gurney,Lin, Mei1}. When $f(u,v)=-au+b(1-u)v$, equation (\ref{eqcon1}) is called the \emph{vector disease model} as proposed by Ruan and Xiao~\cite{Ruan}. When $f(u,v)=bv\exp\{-\gamma\tau\}-\delta u^2$ and $h(u)=\frac{1}{\sqrt{4\pi     \alpha\tau}}\exp\{\frac{-y^2}{4\alpha\tau}\}$, equation (\ref{eqcon1}) is the age-structured reaction diffusion model of a single species proposed by Al-Omari \& Gourley~\cite{Al}. Existence and stability of traveling wavefronts for the reaction-diffusion equation (\ref{eqcon1}) and its special forms has been extensively studied in the literature.

In particular, if $h(x)=\sum^{\infty}_{j=-\infty}J(j)\delta(j)$ with $\delta(\cdot)$ being the Dirac delta function, then (\ref{eq}) become
\begin{equation}
  \label{eq-lde}
  u_t(x,t)= d\cdot\Delta_1u(x,t)+f
  \left(u(x,t),\sum^{\infty}_{j=-\infty}J(j)u(x-j,t-\tau)\right).
\end{equation}
For any fixed $x\in \mathbb{R}$, let $v_n(t)=u(x+n,t)$. Then (\ref{eq-lde}) can be written as the following lattice differential equations
\begin{equation}
  \label{lattice}
  v'_n(t)= d[v_{n+1}(t)+v_{n+1}(t)-2v_{n}(t)]+f
  \left(v_{n}(t),\sum^{\infty}_{j=-\infty}J(j)v_{n-j}(t-\tau)\right),
\end{equation}
where $t>0$ and $n\in \mathbb{Z}$ can be interpreted as particle index. Models involving lattice differential equations can be found in many scientific disciplines, including physical applications~\cite{sco}, chemical reaction theory \cite{Laplante-Erneux}, biology \cite{Bell}, material science \cite{Cahn}, image processing and pattern recognition \cite{Chua}. Moreover, system (\ref{lattice}) and its special cases has been dealt with by Chow \emph{et al.} \cite{Chow}, Ma and Zou \cite{Ma}, Wu and Zou \cite{WuZou2}.
From an earlier work of Keener \cite{Keener}, one knows that as far as traveling wave fronts
are concerned, a discrete model could behave totally different from its continuous version. It is such an essential difference that drives us to investigate the
existence, uniqueness, and asymptotic stability of traveling wave fronts of (\ref{eq}).

We are interested in wave propagation phenomena, more specifically in monotone traveling waves $u(x, t) = \phi(x + ct)$ for (\ref{eq}), with $\phi$ saturating at 0 and $K$. As usual, we call $c$ the \emph{traveling wave speed} and $\phi$ the \emph{profile} of the wavefront. Then a traveling wave $\phi(\xi)$, where $\xi= x + ct$, is a solution of the associated traveling wave equation
\begin{equation}
  \label{eq3}
  \left\{\begin{split}
      -c\phi'(\xi)+d\cdot \Delta_1\phi(\xi)+f(\phi(\xi),
      (h*\phi)(\xi-c\tau))=0,\quad \xi\in \mathbb{R},\\
      \lim\limits_{\xi\to-\infty}\phi(\xi)=0,\qquad
      \lim\limits_{\xi\to \infty}\phi(\xi)=K,\\
      0\leq \phi(\xi)\leq K,\quad \xi\in \mathbb{R},
    \end{split}\right.
\end{equation}
where $(h*\phi)(\xi)=\int_{\mathbb{R}}h(y)\phi(\xi-y)dy$. In contrast to the second order equation~\eqref{eq4}, the differential-difference equation~\eqref{eq3} is a genuinely infinite-dimensional problem.  As described in the work of Rustichini~\cite{Rustichini1,Rustichini2}, the natural setting for an equation such as~\eqref{eq3} is the Banach space $C[-1, 1]$ of continuous functions on $[-1, 1]$.  While much is known about retarded (delay) differential equations (such as involving $\varphi(\xi-1)$ but not $\varphi(\xi + 1)$), much less is known about so-called \emph{mixed}-type equations as~\eqref{eq3}, where both the forward shift $\xi+ 1$ and the backward shift $\xi-1$ appear. This difference makes the study of~\eqref{eq3} more different and more complicated. This is in particular true in comparison to traveling waves in partial differential equations. For example, even if the Laplacian is augmented by difference-differential terms of the kind considered here, the technical difficulties in proving the existence of travelling waves increase considerably~\cite{Hupkes09a}. An example of the difficulties one faces when studying forward-backward equations is that given an initial condition (here on $[-1,1]$), there is in general no solution on the real line, and sophisticated Fredholm arguments are needed~\cite{Harterich:02a}.

We mention two special cases covered by the following analysis: If $f(u,v)=-au+b(v)$ and $h(x)=\delta(x)$, where $b\in C^1([0,\infty],\mathbb{R})$ and $\delta(\cdot)$ is the Dirac delta function, then (\ref{eq}) reduces to the local equation
\begin{equation}
  \label{ex0}
  u_t(x,t)=d\cdot\Delta_1u(x,t)-au(x,t)+b(u(x,t-\tau)),\quad x\in
  \mathbb{R},\,t\geq 0.
\end{equation}
If $f(u,v)=g(u)$ and $g(u)$ denotes a Lipschitz continuous function satisfying $g(u)>0=g(0)=g(1)$ for all $u\in (0,1)$, then (\ref{eq}) becomes
\begin{equation}
  \label{eqc2}
  u_t(x,t)= d\cdot\Delta_1u(x,t)+g(u(x,t)).
\end{equation}

Before recalling the precise result and recalling relevant existence and stability theorems, we state the structural assumptions required for the stability result. Namely, we require the nonlinear functions $f(u,v)$ and $h(u)$ to satisfy the hypotheses
\begin{description}
\item[(F1)] $f\in C([0,K]^2,\mathbb{R})$, $f(0,0)=f(K,K)=0$,   $f(u,u)>0$ for all $u\in (0,K)$, $\partial_2f(u,v)\geq 0$ for all   $(u,v)\in [0,K]^2$, where $K$ is a positive constant.
\item[(F2)] There exist some $M>0$ and $\sigma\in (0,1]$ such that
   \begin{equation*}
     0\leq \partial_1f(0,0)u+\partial_2f(0,0)v-f(u,v)\leq              M(u+v)^{1+\sigma}
   \end{equation*}
   for all $(u,v)\in [0,K]^2$ and    $\partial_1f(K,K)+\partial_2f(K,K)<0$.
\item[(F3)] $\partial_{1}f(0,0)\leq 0$ and $\partial_{ij}f(u,v)\leq 0$   ($i,j=1,2$) for all $(u,v)\in [0,K]^2$.
 \item[(H1)] $h (x)$ is nonnegative, even, integrable and satisfies    $\int_{\mathbb{R}}h(x)dx=1$.
 \item[(H2)] There exists some $\lambda_0>0$ (possibly equal to    $\infty$) such that $\int^{\infty}_0h(x)\exp\{\lambda x\}dx<\infty$    for all $\lambda<\lambda_0$.
\end{description}

We remark that assumptions (F1) and (F2) are standard; (F1) shows that (\ref{eq}) has two equilibria $0$ and $K$. Furthermore, condition (F2) together with (F1) implies that $\partial_1f(0,0)+\partial_2f(0,0)\geq \frac{2}{K}f(\frac{K}{2},\frac{K}{2})>0$, hence 0 is unstable and $K$ is stable.

We point to some related literature, giving a very incomplete synopsis. Regarding work on related equations, we only mention the analysis of equations with a local nonlinearity, but general nonlocal expressions instead of the spatial discrete Laplacian (e.g.,~\cite{Coville}), and the interesting work on a Fisher-KPP equation with a non-local saturation effect, where no maximum principle holds~\cite{Berestycky}. As for stability, Sattinger \cite{Sattinger} used the spectrum-analysis method to prove the wave stability for the Fisher-KPP nonlinearity, when the initial perturbation has an exponential decay. Schaaf \cite{Schaaf} proved linearized stability for the Fisher-KPP nonlinearity by a spectral method.  Smith and Zhao \cite{SmithZhao} considered a \emph{bi-stable} nonlinearity, and proved a global stability result for the traveling wave solution. Stability of wavefronts of reaction-diffusion equations with critical speeds was studied by Kirchg\"assner \cite{Kirchgassner} and Gallay \cite{Gallay}, Mei and his coworkers \cite{Mei1,Mei2}.

Concerning the stability of traveling waves of the model (\ref{eq}), an effective method is the comparison principle combined with the squeezing technique, which has been used by many authors for various monostable equation (see, for example, \cite{Chen1997,LvWang2010,Ma,Wang}).  Another effective method is the (technical) weighted energy method.  This method was used by Mei \emph{et al.} \cite{Mei2} for the Nicholson's blowflies equation, i.e., equation (\ref{eqcon1}) with $f(u,v)=-\gamma u+pv\exp\{-av\}$ and $h(u)=\delta(u)$, and further employed by many researchers to prove the stability of monotone traveling waves of various monostable reaction-diffusion equations with delays. See, e.g., \cite{Gourley2005,Huang,LiMei,LvWang2012,Mei3} and the references therein. For this method, the key step is to establish \emph{a priori} estimates. It is natural to ask if the method can be extended to spatially discrete reaction-diffusion equations with nonlocal delay effects. We give here an affirmative answer.
To the best of our knowledge, this is the first time the asymptotic stability of monotone traveling waves of discrete reaction-diffusion equations with nonlocal delay effects has been studied. The proofs rely on the weighted energy method.

\subsection{Existence and asymptotic stability}
\label{sec:Exist-asympt-stab}
The existence of wave fronts and its spreading speed can be selected linearly or nonlinearly; see, e.g.,~\cite{Lucia}.  The main idea of Lucia \emph{et al.}~\cite{Lucia} is that the linearized system near the unstable equilibrium and its corresponding characteristic equation define a (small) speed $c=c_0$. For a monotone system, Liang and Zhao~\cite{Liang1,Liang2} proved the existence of traveling waves of reaction-diffusion equations for $c\geq c^{*}$.  It is known that $c^{*}\geq c_0$.  We say \emph{linear selection} holds if $c^{*}=c_0$ and \emph{nonlinear selection} takes place if $c^{*}>c_0$. Usually, to have linear selection, the nonlinearity $f$ in the model is required to be sublinear so that the critical speed $c^{*}=c_0$ can be determined by the characteristic equation. In this paper, condition (F2) implies that the nonlinearity $f$ is sublinear. Thus, our model has a linear selection.

In a recent paper~\cite{Guo:11a}, we have shown the existence and uniqueness of traveling waves for (\ref{eq}). We recall the main result of~\cite{Guo:11a} and point out that (F3) is not required for the results in this subsection.

\begin{theorem}\label{thm-existence} Under assumptions (F1), (F2), (H1), and (H2), there exists a minimal wave speed $c^*>0$ such that   for each $c\geq c^*$, equation (\ref{eq}) has a traveling wavefront   $\phi(x+ct)$ satisfying (\ref{eq3}). Moreover,
  \begin{enumerate}
  \item the solution $\phi$ of (\ref{eq3}) is unique up to a     translation.
  \item Every solution $\phi$ of (\ref{eq3}) is strictly monotone,     i.e., $\phi'(\xi)>0$ for all $\xi\in \mathbb{R}$.
  \item Every solution $\phi$ of (\ref{eq3}) satisfies $0 <\phi(\cdot)     <K$ on $\mathbb{R}$.
  \item Any solution of (\ref{eq3}) satisfies     $\lim_{\xi\to-\infty}\phi'(\xi)/\phi(\xi)=\lambda$, with $\lambda$     being the minimal positive root of
    \begin{equation}\label{1eq}
      c\lambda-d[e^{\lambda}+e^{-\lambda}-2]-\partial_1f(0,0)
      -\partial_2f(0,0)\int_{\mathbb{R}}h(y)
      \exp\{-\lambda             (y+c\tau)\}dy=0.
    \end{equation}
  \item Any solution of (\ref{eq3}) satisfies $\lim_{\xi\to       \infty}\phi'(\xi)/[K-\phi(\xi)]=\gamma$, with $\gamma$ being the     unique root of
    \begin{equation}\label{2eq}
      c\gamma+d[e^{\gamma}+e^{-\gamma}-2]+\partial_1f(K,K)
      +\partial_2f(K,K)\int_{\mathbb{R}}h(y)e^{\gamma(y+c\tau)}dy=0.
    \end{equation}
  \end{enumerate}
\end{theorem}

\begin{theorem}\label{thm-as}
  Under assumptions (F1), (F2), (H1) and (H2), for each solution of   $(c, \varphi)$ of (\ref{eq3}) there exists $\eta=\eta(\varphi)$ such   that   $\lim\limits_{\xi\to-\infty}\frac{\varphi(\xi+\eta)}{e^{\lambda_1(c)\xi}}=1$   for $c>c^*$ and $\lim\limits_{\xi\to-\infty}   \frac{\varphi(\xi+\eta)}{\xi e^{\lambda_1(c)\xi}} =1$ for $c=c^*$.   Moreover,   $\lim\limits_{\xi\to-\infty}\frac{\varphi'(\xi)}{\varphi(\xi)}   =\lambda_1(c)$ for $c\geq c^*$.
\end{theorem}

\begin{theorem}\label{thm-as2}
  Under assumptions (F1), (F2), (H1) and (H2), for each solution of   $(c, \varphi)$ of (\ref{eq3}) there exists $\eta=\eta(\varphi)$ such   that   $\lim\limits_{\xi\to\infty}\frac{K-\varphi(\xi+\eta)}{e^{-\upsilon(c)\xi}}=1$,   where $\upsilon(c)$ is the unique positive zero of   $\widetilde{\Delta}(c,\cdot)$, defined in Lemma~\ref{lem-de2} below.   Moreover,   $\lim\limits_{\xi\to\infty}\frac{\varphi'(\xi)}{K-\varphi(\xi)}=\upsilon(c)$.
\end{theorem}

\subsection{The main statement}
\label{sec:main-statement}

We begin this section with notation and some auxiliary statements. We denote by $L^p(I)$ ($p\geq 1$) the Lebesgue space of integrable functions defined on $I$, and $W^{k,p}(I)$ ($k\geq 0$, $p\geq 1$) is the Sobolev space of the $L^p$-functions $f(x)$ defined on the interval $I$ whose derivatives $\frac{\mathrm{d}^n}{\mathrm{d}x^n}f$ ($n=1,\ldots,k$) also belong to $L^p(I)$. We write $H^k(I)$ for $W^{k,2}(I)$. Let $L^p_w(I )$ ($p\geq 1$) be the weighted $L^p$-space with a weight function $w(x) > 0$ and the norm
\begin{equation*}
  \|f\|_{L^p_w(I)}=\left[\int_{I}w(x)|f(x)|^pdx\right]^{\frac{1}{p}}.
\end{equation*}
Analogously, let $H^k_w(I)$ be the weighted Sobolev space with the norm given by
\begin{equation*}
  \|f\|_{H^k_w(I)}=\left[\sum\limits^k_{n=0}\int_{I}w(x)
    \left|\frac{\mathrm{d}^nf(x)}{\mathrm{d}x^n}\right|^2dx
  \right]^{\frac{1}{2}}.
\end{equation*}

For a given traveling wave $\phi$ of (\ref{eq}) satisfying (\ref{eq3}), define
\begin{equation}\label{ch517}
  \begin{split}
    G_j(\xi)&=\partial_jf\left(\phi(\xi),\int_{\mathbb{R}}
      h(y)\phi(\xi-y-c\tau)dy\right),\quad j=1,2,\\
    B(\xi)&=\int_{\mathbb{R}}h(y)G_2(\xi+y+c\tau)dy.
  \end{split}
\end{equation}
Obviously, $B(\xi)$ and $G_j(\xi)$, $j=1,2$ are non-increasing and satisfy
\begin{equation}\label{ch517*}
  G_1(\infty)=\partial_1f(K,K),\quad
  B(\infty)=G_2(\infty)=\partial_2f(K,K).
\end{equation}

We now introduce some central notations. For $\lambda\in \mathbb{C}$ with $\mathrm{Re}\lambda<\lambda_0$, define a function
\begin{equation*}
  G(\lambda)=\int_{\mathbb{R}}h(y)e^{-\lambda y}dy.
\end{equation*}
Obviously, $G(\lambda)$ is twice differentiable in $[0, \lambda_0)$, moreover, $G(0)=1$, $G'(\lambda)=\int^{\infty}_0yh(y)(e^{\lambda   y}-e^{-\lambda y})dy>0$, and $G''(\lambda)=\int^{\infty}_0y^2h(y)(e^{\lambda y}+e^{-\lambda   y})dy>0$. Set
\begin{equation}\label{eqde0}
  \Delta(c,\lambda)=c\lambda-d[e^{\lambda}+e^{-\lambda}-2]-\partial_1f(0,0)
  -\partial_2f(0,0)e^{-\lambda c\tau}G(\lambda)
\end{equation}
and
\begin{equation}\label{eqde}
  \widetilde{\Delta}(c,\lambda):=c\lambda+d[e^{\lambda}+e^{-\lambda}-2]
  +\partial_1f(K,K)+\partial_2f(K,K)
  e^{\lambda         c\tau}G(-\lambda)=0
\end{equation}
for all $c\in \mathbb{R}$ and $\lambda\in \mathbb{C}$ with $c\geq 0$ and $\mathrm{Re}\lambda<\lambda^+$, where $\lambda^+=\lambda_0$ if $\partial_2f(0,0)>0$ and $\lambda^+=+\infty$ if $\partial_2f(0,0)=0$.

We recall two elementary technical statements from \cite{Guo:11a}.

\begin{lemma}\label{lem2}
  There exist $c^*>0$ and $\lambda^*\in (0,\lambda^+)$ such that   $\Delta(c^*,\lambda^*)=0$ and $\Delta_{\lambda}(c^*,\lambda^*)=0$.   Furthermore,
  \begin{description}
  \item[(i)] if $0<c<c^*$, then $\Delta(c,\lambda)<0$ for all     $\lambda\geq 0$;
  \item[(ii)] if $c>c^*$, then the equation $\Delta(c,\cdot)=0$ has     two positive real roots $\lambda_1(c)$ and $\lambda_2(c)$ with     $0<\lambda_1(c)<\lambda^*<\lambda_2(c)<\lambda^+$ such that     $\lambda'_1(c)<0$, $\lambda'_2(c)>0$, $\Delta(c,\lambda)>0$ for     all $\lambda\in (\lambda_1(c),\lambda_2(c))$, and     $\Delta(c,\lambda)<0$ for all $(-\infty,\lambda^+)\setminus     [\lambda_1(c),\lambda_2(c)]$.
  \end{description}
\end{lemma}

\begin{lemma}\label{lem-de2} Under assumptions (F1) and (F2), for each     fixed $c\geq 0$, $\widetilde{\Delta}(c,\cdot)$ has exactly one   positive zero $\upsilon(c)$.
\end{lemma}

With the definitions made in these lemmas, we can introduce
\begin{equation*}
  \mathcal{M}(c,\mu)=c\lambda-d(e^{\lambda}+e^{-\lambda}-2)
  -\partial_1f(0,0)-\mu-\partial_2f(0,0)e^{(\mu-c\lambda)\tau}G(\lambda)
\end{equation*}
and
\begin{equation*}
  \mathcal{N}(\mu)= \mu+\partial_1f(K,K)+e^{\mu\tau}\partial_2f(K,K),
\end{equation*}
where $\lambda$ is any fixed number in $(\lambda_1(c), \lambda^*)$ if $c>c^*$ and $\lambda=\lambda^*$ if $c=c^*$. We notice $\mathcal{M}(c,0)=\Delta(c,\lambda)$, which is positive if $c>c^*$ and equal to 0 if $c=c^*$. Thus, if $c>c^*$, there exists $\mu_0>0$ such that $\mathcal{M}(c,\mu)>0$ for all $\mu\in [0,\mu_0)$. In view of $\partial_1f(K,K)+\partial_2f(K,K)<0$ (by assumption (F2)), there exists a small $\mu>0$ so that $\mathcal{N}(\mu)<0$.

In view of (\ref{ch517*}), there exists a sufficiently large $x_0$ such that
 \begin{equation}\label{assume2}
   \mu+G_1(\xi)+e^{\mu\tau}B(\xi)<0\quad
   \mbox{for all $\xi\in      [x_0,\infty)$}.
 \end{equation}
 For this choice of $x_0$ and any given $c\geq c^*$, we define
\begin{equation}\label{weight}
  w(x)=\begin{cases}
      e^{-\lambda(x-x_0)} &\mbox{for $x\leq x_0$},\\
      1 &\mbox{for $x> x_0$},
    \end{cases}
\end{equation}
where $\lambda$ is any fixed number in $(\lambda_1(c), \lambda^*)$ if $c>c^*$ and $\lambda=\lambda^*$ if $c=c^*$. It is easy to see that $w(x)\geq 1$ for all $x\in \mathbb{R}$ and $w(-\infty)=\infty$. We will use this $w$ for the weighted Lebesgue and Sobolev spaces.

We now state the main result of this article.
\begin{theorem}\label{thm-sta}
  Under assumptions (F1)--(F3) and (H1)--(H2), for a given traveling   wave $\phi$ of (\ref{eq}) satisfying (\ref{eq3}), if the initial   data satisfies
  \begin{equation*}
    \mbox{$0\leq u_0(x,s)\leq K$ for all $(x,s)\in
      \mathbb{R}\times       [-\tau,0]$},
  \end{equation*}
  and the initial perturbation $u_0(x,s)-\phi(x+cs)$ is in   $C([-\tau,0], H^1_w(\mathbb{R})\cap H^1(\mathbb{R}))$, then the   solution of (\ref{eq}) with initial date $u_0(x,s)$, $(x,s)\in   \mathbb{R}\times [-\tau,0]$, uniquely exists and satisfies
  \begin{equation*}
    \mbox{$0\leq u(x,t)\leq K$ for all $(x,t)\in
      \mathbb{R}\times       [0,\infty)$}
  \end{equation*}
  and
  \begin{equation*}
    u(x,t)-\phi(x+ct)\in C([0,\infty),
    H^1_w(\mathbb{R})\cap     H^1(\mathbb{R})).
  \end{equation*}
  Moreover, when $c>c^*$, the solution $u(x,t)$ converges to the   noncritical traveling wave $\phi(x+ct)$ exponentially,
  \begin{equation*}
    \sup\limits_{x\in \mathbb{R}}|u(x,t)-\phi(x+ct)
    |\leq Ce^{-\mu t},\quad t>0,
  \end{equation*}
  for a positive constant $\mu$, where $\mu=\mu(c,\lambda)$ for   $\lambda\in (\lambda_1(c),\lambda^*]$ satisfies   $\mathcal{M}(c,3\mu)>0$ and $\mathcal{N}(3\mu)<0$.  When $c=c^*$,   the solution $u(t, x)$ converges to the critical traveling wave $\phi(x + ct)$ algebraically,
  \begin{equation*}
    \sup\limits_{x\in \mathbb{R}}|u(x,t)
    -\phi(x+ct)|\leq     C(1+t)^{-1/2},\quad t>0.
  \end{equation*}
\end{theorem}
This theorem is proved in Sections \ref{sec:Global-Existence} and~\ref{sec:Global-Stability}.

We remark that the stability of the critical traveling wave solutions to either local or nonlocal time-delayed discrete reaction-diffusion equations is a challenging problem. It is well known that the stability of the critical waves is very important in the study of biological invasions (see, e.g., \cite{Thieme,Liang1,Liang2}). By a careful inspection of the traveling wave equation, using the concavity of the function $f(\cdot,\cdot)$, we first establish a weighted $L^1$- and $L^2$-energy estimate of solutions and then obtain the desired $L^2$-energy estimate as well as the exponential convergence rate to the noncritical wave by the ordinary weighted energy method. When the wave is critical, the convergence rate to the wave is proved to be algebraic for the equivalent integral equation with a time-delayed Green function. This seems to be the first result on stability of the critical traveling wave solutions to time-delayed discrete reaction-diffusion equations.

We briefly comment on key differences to the case where the classical Laplacian is investigated instead of the spatially discrete Laplacian. Obviously, the Greens functions arguments differ. Also, the assumptions on the nonlinearity made in~\cite{Mei4,Mei4a} differ from the ones we make; therefore the key energy estimates differ, as discussed in Remark~\ref{diff-Mei}. In particular, the weaker assumption we make on $f$ leads to the estimate in Lemma~\ref{lem-ch51-0}.

Theorem \ref{thm-sta} also shows how the time-delay $\tau$ effects the convergence rate $\mu$ to the non-critical traveling waves. The effect of the time-delay will essentially make the decay rate $\mu$ of the solutions slow down. Namely, $\mu$ becomes the biggest as $\tau\to 0$ and $\mu$ tends to the smallest $0$ as $\tau\to\infty$. These stability results improve and develop the existing works on monostable waves.

This paper is organized as follows. In Section \ref{sec:Global-Existence}, we prove the global existence and uniqueness of solutions to (\ref{eq}). Section~\ref{sec:Global-Stability} is devoted to proving that the traveling wavefronts are exponentially stable with respect to perturbations in some exponentially weighted $L^{\infty}$ spaces, and we obtain the time decay rates of $\sup_{x\in\mathbb{R}}|u(x, t)-\phi(x + ct)|$ by a weighted energy estimate. As the applications of this main result, in Section~\ref{sec:Applications}, we obtain the global and exponential stability of all noncritical traveling waves and the algebraic stability of the critical wave for a host-vector disease model, an age-structured population model, and a nonlocal Nicholson's blowflies model.

\section{Global existence}
\label{sec:Global-Existence}

This section is devoted to the proofs of the global existence and uniqueness of a solution to (\ref{eq}). Throughout this section, we assume that (F1)--(F3) and (H1)--(H2) hold. We consider the initial value problem
\begin{equation}
  \label{eqs41}
  \left\{
    \begin{array}{rcl}
      u_t(x,t) &=& F[u](x,t),\quad x\in \mathbb{R},\, t>0,\\
      u(t,s) &=& \varphi(x,s),\quad x\in \mathbb{R},\, s\in [-\tau,0],
    \end{array}\right.
\end{equation}
where
\begin{equation*}
  F[u](x,t)=d\cdot\Delta_1u(x,t)+f\left(u(x,t),(h*u)(x,t-\tau)\right).
\end{equation*}
Let $X = BUC (\mathbb{R},\mathbb{R})$ be the Banach space of all bounded and uniformly continuous functions from $\mathbb{R}$ to $\mathbb{R}$ equipped with the supremum norm $\|\cdot\|_X$. Clearly, (\ref{eqs41}) is equivalent to
\begin{equation}
  \label{eqs42}
    u(x,t) = \varphi(x,0)e^{-\mu t}+\int^t_0e^{\mu(s-t)}\{\mu
    u(x,s)+F[u](x,s)\}ds,
\end{equation}
where $\mu>0$ is introduced to ensure the monotonicity of $\mu u+F[u]$ in $u$.  The existence of solutions then follows by Picard's iteration (see, e.g., \cite{Hale}).  To examine stability, we make use of the following comparison result.

\begin{lemma}
  \label{CP}
  Assume that $u^1$ and $u^2$ are continuous functions on   $\mathbb{R}\times [0, \infty)$ such that $0\leq u^1, u^2\leq K$ on   $\mathbb{R}\times [0, \infty)$, that $u^1\geq u^2$ on   $\mathbb{R}\times [-\tau,0]$, and that
  \begin{equation}
    \label{eqs44}
    u^1_t(x,t)-F[u^1](x,t)\geq u^2_t(x,t)-F[u^2](x,t)
  \end{equation}
  for all $x\in \mathbb{R}$ and $t>0$. Then $u^1\geq u^2$ on   $\mathbb{R}\times (0,\infty)$.
\end{lemma}

\proof Let $\mu$ be such that $\mu-\partial_1f(0,0)-\partial_2 f(0,0)e^{-\mu\tau}>0$. Since $w(x,t)=u^2(x,t)-u^1(x,t)$ is continuous and bounded, $\omega(t)=\sup_{\mathbb{R}}w(\cdot,t)$ is continuous on $[0,\infty)$. Suppose the assertion is not true. Then there exists $t_0>0$ such that $\omega(t_0)>0$; it is no loss of generality to assume that
\begin{equation*}
  \omega(t_0)e^{-\mu t_0} >\omega(s)e^{-\mu s}\quad
  \mbox{for all     $s\in [-\tau,t_0)$}.
\end{equation*}
Let $\{x_j\}^{\infty}_{j=1}$ be a sequence on $\mathbb{R}$ such that $w(x_j, t_0)>0$ for all $j\geq 1$ and $\lim\limits_{j\to   \infty}w(x_j,t_0)=\omega(t_0)$. Let $\{t_j\}^{\infty}_{j=1}$ be a sequence in $[0,t_0]$ such that
\begin{equation*}
  e^{-\mu t_j}w(x_j,t_j)=\max\limits_{t\in [0,t_0]}
  \left\{e^{-\mu       t}w(x_j,t)\right\}.
\end{equation*}
Since $\omega (t)e^{-\mu t} <\omega (t_0)e^{-\mu t_0}$ for all $t\in [0,t_0)$, we have $\lim\limits_{j\to\infty}t_j=t_0$ and $\lim\limits_{j\to\infty}w(x_j,t_j)=\omega(t_0)$. In addition, for each $j\geq 1$, we obtain
\begin{equation*}
  0\leq \left.\underline{D}_t\left\{e^{-\mu         t}
      w(x_j,t)\right\}\right|_{t=t_j^-}=e^{-\mu t_j}
  \left\{\underline{D}_tw(x_j,t_j)-\mu w(x_j,t_j)\right\},
\end{equation*}
where $\underline{D}_tu(x,t)=\liminf\limits_{h\to   0}\frac{u(x,t)-u(x,t-h)}{h}$.  Hence, $\underline{D}_t(u^2-u^1)(x_j,t_j)=\underline{D}_tw(x_j,t_j)\geq \mu w(x_j,t_j)$. Thus, it follows from (\ref{eqs44}) and (F3) that
\begin{align*}
    0 &\geq  \underline{D}_tw(x_j,t_j)-d\cdot\Delta_1w(x_j,t_j)\\
    &\quad{}+f(u^1(x_j,t_j),h*u^1(x_j,t_j-\tau))-f(u^2(x_j,t_j),
    h*u^2(x_j,t_j-\tau))  \\
    &\geq  (\mu +2d)w(x_j,t_j)-d[w(x_j-1,t_j)+w(x_j+1,t_j)]\\
    &\quad{}-\partial_1f(0,0)w(x_j,t_j)-\partial_2f(0,0)h*w(x_j,t_j-\tau)\\
    &\geq  [\mu     +2d-\partial_1f(0,0)]
    w(x_j,t_j)-2d\omega(t_j)-\partial_2f(0,0)\omega(t_j-\tau).
\end{align*}
Letting $j\to \infty$, we have
\begin{align*}
    0 &\geq [\mu     +2d-\partial_1f(0,0)]\omega(t_0)
    -2d\omega(t_0)-\partial_2f(0,0)\omega(t_0-\tau)\\
    &\geq  [\mu -\partial_1f(0,0)-\partial_2f(0,0)
    e^{-\mu         \tau}]\omega(t_0).
\end{align*}
Recalling that $\mu -\partial_1f(0,0)-\partial_2f(0,0)e^{-\mu   \tau}>0$, we conclude that $\omega(t_0)\leq 0$, which contradicts $\omega(t_0)>0$. Therefore, $w(x,t)=u^2(x,t)-u^1(x,t)\geq 0$ for all $(x,t)\in \mathbb{R}\times [0,\infty)$. This completes the proof. \qed

\begin{theorem}\label{thm-sta2}
  For a given traveling   wave $\phi$ of (\ref{eq}) satisfying (\ref{eq3}), if the initial   data satisfies $0\leq u_0(x,s)\leq K$ for all $(x,s)\in
      \mathbb{R}\times       [-\tau,0]$,
  and the initial perturbation $u_0(x,s)-\phi(x+cs)$ is in   $C([-\tau,0], H^1(\mathbb{R}))$, then the   solution of (\ref{eq}) with initial date $u_0(x,s)$, $(x,s)\in   \mathbb{R}\times [-\tau,0]$, uniquely exists and satisfies
  $0\leq u(x,t)\leq K$ for all $(x,t)\in
      \mathbb{R}\times       [0,\infty)$
  and
  $u(x,t)-\phi(x+ct)\in C([0,\infty),
    H^1(\mathbb{R}))$.
\end{theorem}

The rest of this section is devoted to the proof of Theorem \ref{thm-sta2}. Firstly,
it follows from the comparison principle (Lemma \ref{CP}) that $0\leq u(x,t)\leq K$
for all $(x,t)\in \mathbb{R}\times [0,\infty)$. Let $v(x,t)=u(x,t)-\phi(x+ct)$. Then (\ref{eq}) can be written as
\begin{equation}\label{eq-in}
\aligned
&v_t(x,t)=d\cdot\Delta_1v(x,t)+\mathcal{Q}(x,t),\quad
&(x,t)\in \mathbb{R}\times [0,\infty),\\
&v(x,s)=v_0(x,s),\quad &(x,s)\in \mathbb{R}\times [-\tau,0],
\endaligned
\end{equation}
where $v_0(x,s)=u(x,s)-\phi(x+cs)$ and
$$
\mathcal{Q}(x,t)=f\left(u(x,t),(h*u)(x,t-\tau)\right)
-f\left(\phi(x+ct),(h*\phi)(x+ct-c\tau)\right).
$$
The following result on the local
existence, uniqueness and extension of solutions is standard. It can be proved using
the standard iteration method (see, e.g., \cite{Hale}). The proof is thus omitted.

\begin{lemma}\label{local-ex} Given
$v_0\in C([-\tau,0], H^1(\mathbb{R}))$,
 there exists $t_0 > 0$ such that problem (\ref{eq-in}) has a unique solution $v \in C([0,t_0), H^1(\mathbb{R}))$. Furthermore, let $[0,T)$ be the maximal interval of existence and
$v \in C([0,T), H^1(\mathbb{R}))$. Then either $T =+\infty$ or the solution blows up in finite time, in the sense that $T<+\infty$ and $\lim_{t\to T^-}\|v(\cdot,t)\|_{H^1(\mathbb{R})}=\infty$.
\end{lemma}

The  lemma above implies the local existence of a solution in the stated solution space for a time interval $[0,T)$, for some $T> 0$. Then, the solution either  exists globally in the given function space or blows up in finite time, in the $H^1$ norm. We now further show, using the energy method, for any time $T> 0$, that the solution defined on $(0, T)$ in the given space is bounded by a constant depending on $T$ and does not blow up. Consequently, we obtain the global existence in the given function space.

\begin{lemma}\label{bounded} In the setting of Lemma~\ref{local-ex}, let $v$ be a solution in $C([0,T), H^1(\mathbb{R}))$ for
$0 <T<\infty$. Then there exists positive constant $C_0$, independent of $T$, such that
$$
\|v(\cdot,t)\|_{H^1(\mathbb{R})}\leq C_0\left(\|v_0(\cdot,0)\|^2_{H^1(\mathbb{R})}
+\int^{0}_{-\tau}\|v_0(\cdot,s)\|^2_{H^1(\mathbb{R})}\mathrm{d}s \right)e^{2\vartheta_1t}
$$
for all $t\in [0,T)$, where $\vartheta_1:=\max\{\left|\partial_jf\left(u,v\right)\right| :\,u,v\in [0,K],\, j=1,2\}$.
\end{lemma}

\proof It follows from assumption (F3) and the mean-value theorem that
$$
|\mathcal{Q}(x,t)|\leq \vartheta_1|v(x,t)|+\vartheta_1|(h*v)(x,t-\tau)|.
$$
Multiplying (\ref{eq-in}) by $v(x,t)$, integrating over $\mathbb{R}\times [0,t]$, $t\in [0,T)$, and noticing that
\begin{equation}\label{ch6-6}
  v(\xi)\Delta_1v(\xi)\leq \frac 1 2 \Delta_1v^2(\xi)
  \quad (\mbox{using the binomial formula}),
\end{equation}
 we have
\begin{align*}
\frac 1 2 \|v(\cdot,t)\|^2_{L^2(\mathbb{R})}
\leq &
\frac 1 2 \|v(\cdot,0)\|^2_{L^2(\mathbb{R})}
+\vartheta_1\int^t_0\|v(\cdot,s)\|^2_{L^2(\mathbb{R})}\mathrm{d}s\\
&+\vartheta_1\int^t_0\int_{\mathbb{R}}
\int_{\mathbb{R}}h(x-y)|v(x,s)v(y,s-\tau)|\mathrm{d}y\mathrm{d}x\mathrm{d}s\\
\leq &
\frac 1 2 \|v(\cdot,0)\|^2_{L^2(\mathbb{R})}
+\vartheta_1\int^t_0\|v(\cdot,s)\|^2_{L^2(\mathbb{R})}\mathrm{d}s\\
&+\frac{1}{2}\vartheta_1\int^t_0\int_{\mathbb{R}}\int_{\mathbb{R}}h(x-y)\left[v^2(x,s)+v^2(y,s-\tau)\right]\mathrm{d}y\mathrm{d}x\mathrm{d}s\\
=&
\frac 1 2 \|v(\cdot,0)\|^2_{L^2(\mathbb{R})}
+\frac{3}{2}\vartheta_1\int^t_0\|v(\cdot,s)\|^2_{L^2(\mathbb{R})}\mathrm{d}s
+\frac{1}{2}\vartheta_1\int^{t-\tau}_{-\tau}\|v(\cdot,s)\|^2_{L^2(\mathbb{R})}\mathrm{d}s\\
\leq &
\frac 1 2 \|v(\cdot,0)\|^2_{L^2(\mathbb{R})}+\frac{1}{2}\vartheta_1\int^{0}_{-\tau}\|v(\cdot,s)\|^2_{L^2(\mathbb{R})}\mathrm{d}s
+2\vartheta_1\int^t_0\|v(\cdot,s)\|^2_{L^2(\mathbb{R})}\mathrm{d}s.
\end{align*}
Applying Gronwall's inequality, we have
$$
\frac 1 2 \|v(\cdot,t)\|^2_{L^2(\mathbb{R})}\leq \left(\frac 1 2 \|v_0(\cdot,0)\|^2_{L^2(\mathbb{R})}+\frac{1}{2}\vartheta_1\int^{0}_{-\tau}\|v_0(\cdot,s)\|^2_{L^2(\mathbb{R})}\mathrm{d}s \right)e^{2\vartheta_1t}
$$
for all $t\in [0,T)$.

Let us differentiate (\ref{eq-in}) with respect to $x$, then we have
\begin{equation}\label{eq-in2}
\aligned
&(v_x)_t(x,t)=d\cdot\Delta_1v_x(x,t)+\mathcal{Q}_x(x,t),\quad
&(x,t)\in \mathbb{R}\times [0,\infty),\\
&v(x,s)=v_0(x,s)\quad &(x,s)\in \mathbb{R}\times [-\tau,0],
\endaligned
\end{equation}
where
\begin{align*}
\mathcal{Q}_x(x,t)
=&\partial_1f\left(u(x,t),(h*u)(x,t-\tau)\right)u_x(x,t)\\
&+\partial_2f\left(u(x,t),(h*u)(x,t-\tau)\right)(h*u_x)(x,t-\tau)\\
&-\partial_1f\left(\phi(x+ct),(h*\phi)(x+ct-c\tau)\right)\phi'(x+ct)\\
&-\partial_2f\left(\phi(x+ct),(h*\phi)(x+ct-c\tau)\right)(h*\phi')(x+ct-c\tau).
\end{align*}
It follows from assumption (F3) and the mean-value theorem that
\begin{align*}
|\mathcal{Q}_x(x,t)|
\leq &\vartheta_1|v_x(x,t)|+\vartheta_1|(h*v_x)(x,t-\tau)|\\
&+\vartheta_2\left[|v(x,t)|+|(h*v)(x,t-\tau)|\right]|\phi'(x+ct)|\\
&+\vartheta_2\left[|v(x,t)|+|(h*v)(x,t-\tau)|\right]|(h*\phi')(x+ct-c\tau)|.
\end{align*}
where $\vartheta_2=\max\{\partial_{jk}f\left(u,v\right):\,u,v\in [0,K],\, j,k=1,2\}$. Similarly, multiplying (\ref{eq-in2}) by $v_{x}(x,t)$, integrating over $\mathbb{R}\times [0,t]$, $t\in [0,T)$, and then applying Gronwall's inequality, we have
$$
\frac 1 2 \|v_x(\cdot,t)\|^2_{L^2(\mathbb{R})}\leq C_1\left(\frac 1 2 \|v_0(\cdot,0)\|^2_{H^1(\mathbb{R})}+\int^{0}_{-\tau}\|v_0(\cdot,s)\|^2_{H^1(\mathbb{R})}\mathrm{d}s \right)e^{2\vartheta_1t}
$$
for some positive constant $C_1 > 0$ and all $t\in [0,T)$.
This completes the proof. \qed

Finally, Theorem \ref{thm-sta2} follows immediately from the previous lemmas.

\section{Global stability}
\label{sec:Global-Stability}

This section is devoted to proving the stability of all noncritical traveling waves to (\ref{eq}) with an exponential convergence rate, and the stability of the critical traveling wave with an algebraic convergence rate. Throughout this section, it is assumed that (F1)--(F3) and (H1)--(H2) hold.  We first give some auxiliary statements for linear delayed ODEs.

\begin{lemma}[\cite{Khusainov,Mei-Wang}]\label{lem-KH}
  Let $z(t)$ be the solution to the linear system with time delay   $\tau>0$,
  \begin{equation}\label{delay}
    \begin{split}
      z'(t)+c_1z(t)=c_2z(t-\tau),\\
      z(s)=z_0(s),\quad s\in [-\tau,0].
    \end{split}
  \end{equation}
  Then
  \begin{equation}\label{delay1}
    z(t)=e^{-c_1(t+\tau)}e^{c_3t}_{\tau}z_0(-\tau)
    +\int^0_{-\tau}e^{-c_1(t-s)}e^{c_3(t-\tau-s)}_{\tau}
    \left[z'_0(s)+c_1z_0(s)\right]ds,
  \end{equation}
  where $c_3=c_2e^{c_1\tau}$ and $e^{c_3t}_{\tau}$ is so-called   delayed exponential function in the form
  \begin{equation*}
    e^{c_3t}_{\tau}=
    \begin{cases}
      0 & -\infty<t<-\tau,\\
      1 & -\tau\leq t<0,\\
      1+\frac{1}{1!}\exp\{c_3t\} 0\leq t<\tau,\\
      1+\frac{1}{1!}\exp\{c_3t\}+\frac{1}{2!}\exp\{c_3^2(t-\tau)^2\}
      & \tau\leq t<2\tau,\\
      \vdots & \vdots\\
      1+\frac{1}{1!}\exp\{c_3t\}+\cdots+
      \frac{1}{m!}\exp\{c_3^m[t-(m-1)\tau]^m\},
      & (m-1)\tau\leq t<m\tau,\\
      \vdots & \vdots
    \end{cases}.
  \end{equation*}
  Furthermore, when $c_1\geq c_2\geq 0$, there exist constants $C>0$   and $\varepsilon\in (0,1)$ such that
  \begin{equation*}
    e^{-c_1t}e^{c_3t}_{\tau}\leq Ce^{-\varepsilon(c_1-c_2)t},
    \quad     t>0.
  \end{equation*}
\end{lemma}

\begin{lemma}\label{lem-G} For all $\xi\in \mathbb{R}$ and $t\geq
  0$,
  \begin{equation}\label{stand2}
    \frac{1}{2\pi}\int_{\mathbb{R}}\exp\{-2td\varepsilon
    \cosh(\lambda^*)(1-\cos\omega)\} d\omega\leq
    \sqrt{\frac{\pi}{dt\varepsilon}}.
  \end{equation}
\end{lemma}

\proof Define
\begin{equation*}
  v(\xi,t)=\frac{1}{2\pi}\int_{\mathbb{R}}\exp\{tq(\lambda^*,\omega)+i\omega
  \xi\}d\omega,
\end{equation*}
where
\begin{equation*}
  q(\lambda^*,\omega)=d\varepsilon\left[\exp\{\lambda^*+i\omega\}
    +\exp\{-\lambda^*-i\omega\}-\exp\{\lambda^*\}-\exp\{-\lambda^*\}\right].
\end{equation*}
Obviously, $v(\xi,t)$ is the solution to
\begin{equation}\label{stand}
  v_t(\xi,t)=d\varepsilon[e^{\lambda^*}v(\xi+1,t)
  +e^{-\lambda^*}v(\xi-1,t)]-d\varepsilon(e^{\lambda^*}+e^{\lambda^*})v(\xi,t).
\end{equation}
with the initial data $v(\xi,0)$ taken to be the Dirac delta $\delta(\xi)$. For any given $\xi\in \mathbb{R}$, let $v_n(t)=v(\xi+n,t)$ then (\ref{stand}) can be written as
\begin{equation*}
  \partial_t   v_n(t)=d\varepsilon[e^{-\lambda^*}v_{n-1}(t)
  +e^{\lambda^*}v_{n+1}(t)]-d\varepsilon[e^{\lambda^*}+e^{-\lambda^*}]v_{n}(t),
\end{equation*}
where $n$ can be interpreted as particle index. The Green function can be computed as in~\cite{Fath} (see also~\cite{Carpio}): by taking the discrete Fourier transform $F[v](\omega,t) = \sum e^{in\omega} v_n(t)$, we obtain
\begin{equation*}
  \partial_t F[v](\omega,t) =q(\lambda^*,\omega)F[v](\omega,t).
\end{equation*}
The Laplace transform in time yields
\begin{equation*}
  L[F[v]](\omega,s) = \frac{F[v](\omega,0)}{s-q(\lambda^*,\omega)}.
\end{equation*}
Inverting the Laplace transform, we get
\begin{equation*}
  F[v](\omega,t) = F[v](\omega,0) \exp\{tq(\lambda^*,\omega)\}.
\end{equation*}
Taking the inverse discrete Fourier transform then yields (with the initial data taken to be the Dirac delta at the origin)
\begin{align*}
  v_n(t) &= \frac{1}{2\pi} \int^{\pi}_{-\pi}
 \exp\{tq(\lambda^*,\omega)+in\omega\}d\omega\\
 &\leq \frac{1}{2\pi}
 \int^{\pi}_{-\pi}\exp\{2td\varepsilon\cosh(\lambda^*)(\cos\omega-1)\}d\omega\\
 &\leq \frac{1}{2\pi}
 \int^{\pi}_{-\pi}\exp\{-2td\varepsilon\cosh(\lambda^*)
 \sin^2\frac{\omega}{2}\}d\omega.
\end{align*}
Noting that $\frac{2\omega}{\pi}\leq \sin \omega\leq \omega$ for all $\omega\in [0,\frac{\pi}{2}]$, we have
\begin{align*}
  v_n(t)  &\leq \frac{1}{2\pi}\int^{\pi}_{-\pi}
 \exp\{-2td\varepsilon\cosh(\lambda^*)\frac{\omega^2}{\pi^2}\}d\omega\\
 &=\frac{1}{2\sqrt{2td\varepsilon\cosh(\lambda^*)}}
 \int^{\sqrt{2td\varepsilon\cosh(\lambda^*)}}_{-\sqrt{2td\varepsilon
     \cosh(\lambda^*)}}\exp\{-\omega^2\}d\omega\\
 &\leq   \frac{1}{2\sqrt{2td\varepsilon\cosh(\lambda^*)}}
 \int_{\mathbb{R}}\exp\{-\omega^2\}d\omega\\
 &\leq  \sqrt{\frac{\pi}{dt\varepsilon}}.
\end{align*}
Here we have used the fact that $\cosh(\lambda^*)>1$. Hence, for any solution $v(\xi,t)$ to (\ref{stand}), we have
\begin{equation*}
  v(\xi,t)\leq \sqrt{\frac{\pi}{dt\varepsilon}}
\end{equation*}
for all $\xi\in \mathbb{R}$ and $t> 0$. Obviously, the left-hand side of inequality (\ref{stand2}) is exactly $v(0,t)$. Therefore, (\ref{stand2}) holds. This completes the proof. \qed

Let $c\geq c^*$ and the initial data $u_0(x,s)$ be such that $0\leq u_0(x,s)\leq K$ for all $(x,s)\in \mathbb{R}\times [-\tau,0]$, and define
\begin{equation*}
  u^+_0(x,s)=\max\{u_0(x,s),\phi(x+cs)\},\quad
  u^-_0(x,s)=\min\{u_0(x,s),\phi(x+cs)\}
\end{equation*}
for all $(x,s)\in \mathbb{R}\times [-\tau,0]$. Let $u^+(x,t)$ and $u^-(x,t)$ be the corresponding solutions of (\ref{eq}) with the initial data $u^+_0(x,s)$ and $u^-_0(x,s)$ with $(x,s)\in \mathbb{R}\times [-\tau,0]$, respectively. It follows from the comparison principle (Lemma \ref{CP}) that
\begin{equation*}
  0\leq u^-(x,t)\leq u(x,t)\leq u^+(x,t)\leq K
\end{equation*}
and
\begin{equation*}
  0\leq u^-(x,t)\leq \phi(x+ct)\leq u^+(x,t)\leq K
\end{equation*}
for all $(x,t)\in \mathbb{R}\times\mathbb{R}_+$. In what follows, we are going to complete the proof for the stability result (Theorem \ref{thm-sta}) in three steps.

\step{Step 1. The convergence of $u^+(x,t)$ to $\phi(x+ct)$.}

We first prove the convergence of $u^+(x,t)$ to $\phi(x+ct)$. Let $\xi=x+ct$ and
\begin{equation*}
  v(\xi,t)=u^+(x,t)-\phi(x+ct),\quad v_0(\xi,t)=u^+_0(x,s)-\phi(x+cs).
\end{equation*}
It follows that
\begin{equation*}
  v(\xi,t)\geq 0,\quad v_0(\xi,s)\geq 0.
\end{equation*}
Moreover, $v(\xi,t)$ satisfies
\begin{equation}\label{ch51}
    v_t+cv_{\xi}-d\cdot\Delta_1v-\partial_1f(0,0)v-\partial_2f(0,0)
    h*v_{\tau}=Q,
\end{equation}
with the initial data
\begin{equation}\label{ch51'}
  v(\xi,s)=v_0(\xi,s),\quad s\in [-\tau,0],
\end{equation}
where $v=v(\xi,t)$, $v_{\tau}=v(\xi-c\tau,t-\tau)$, $\phi_{\tau}=\phi(\xi-c\tau)$, and
\begin{equation*}
    Q =f\left(v+\phi,h*(v_{\tau}+\phi_{\tau})\right)
    -f\left(\phi,h*\phi_{\tau}\right)-\partial_1f(0,0)v
    -\partial_2f(0,0)h*v_{\tau}.
\end{equation*}
By the Mean Value Theorem, we have
\begin{equation*}
  \begin{split}
    Q =&f\left(v+\phi,h*(v_{\tau}+\phi_{\tau})\right)-
    f\left(\phi,h*(v_{\tau}+\phi_{\tau})\right)\\
    &+f\left(\phi,h*(v_{\tau}+\phi_{\tau})\right)
    -f\left(\phi,h*\phi_{\tau}\right)\\
    &-\partial_1f(0,0)v
    -\partial_2f(0,0)h*v_{\tau}\\
    =&\partial_1f\left(v+\varepsilon_1\phi,h*(v_{\tau}+\phi_{\tau})\right)v-\partial_1f(0,0)v\\
    &+\partial_2f\left(\phi,h*(v_{\tau}+\varepsilon_2\phi_{\tau})\right)h*v_{\tau}-\partial_2f(0,0)h*v_{\tau}\\
    =& \left[\partial_1f\left(v+\varepsilon_1\phi,h*(v_{\tau}+\phi_{\tau})\right)
    -\partial_1f\left(0,h*(v_{\tau}+\phi_{\tau})\right) \right]v\\
    & +\left[\partial_1f\left(0,h*(v_{\tau}+\phi_{\tau})\right)
    -\partial_1f(0,0) \right]v\\
    &+\left[\partial_2f\left(\phi,h*(v_{\tau}+\varepsilon_2\phi_{\tau})\right)
    -\partial_2f\left(0,h*(v_{\tau}+\varepsilon_2\phi_{\tau})\right)\right]h*v_{\tau}\\
        &+\left[\partial_2f\left(0,h*(v_{\tau}+\varepsilon_2\phi_{\tau})\right)
    -\partial_2f\left(0,0\right)\right]h*v_{\tau}\\
    =&\partial_{11}f\left(v+\varepsilon_1\varepsilon_3\phi,h*(v_{\tau}+\phi_{\tau})\right)(v+\varepsilon_1\phi)v\\
    &+\partial_{12}f\left(0,h*(v_{\tau}+\varepsilon_4\phi_{\tau})\right)v[h*(v_{\tau}+\phi_{\tau})]\\
    &+\partial_{21}f\left(\varepsilon_5\phi,h*(v_{\tau}+\varepsilon_2\phi_{\tau})\right)\phi £¨h*v_{\tau}£©\\
    &+\partial_{22}f\left(0,h*(v_{\tau}+\varepsilon_2\varepsilon_6\phi_{\tau})\right)(h*v_{\tau})[h*(v_{\tau}+\varepsilon_2\phi_{\tau})],
  \end{split}
\end{equation*}
where $\varepsilon_1$, $\varepsilon_2$, ..., $\varepsilon_6\in (0,1)$. It follows from assumption (F3) and the nonnegativity of $h$, $\phi$ and $v$ that $Q(\xi,t)\leq 0$.

\begin{lemma}\label{lem-ch51}
  It holds that
  \begin{equation}\label{ch52}
    \|v(t)\|_{L^1_{w}(\mathbb{R})}+\int^{t}_0
    e^{-\mu_1(t-s)}\|v(s)\|_{L^1_{w}(\mathbb{R})}ds
    \leq
    Ce^{-\mu_1t}
  \end{equation}
  for $c>c^*$ and
  \begin{equation}\label{ch52'}
    \|v(t)\|_{L^1_{w}(\mathbb{R})}\leq
    C
  \end{equation}
  for $c=c^*$, where $w(\xi)=e^{-\lambda(\xi-x_0)}$, $\lambda$ is any   fixed number in $(\lambda_1(c), \lambda^*)$ if $c>c^*$ and   $\lambda=\lambda^*$ if $c=c^*$, and $\mu_1>0$ is a constant such   that $\mathcal{M}(c,\mu_1)>0$.
\end{lemma}

\proof In order to establish the energy estimate (\ref{ch52}) and (\ref{ch52'}), sufficient regularity of the solution to (\ref{ch51}) and (\ref{ch51'}) is required. We thus mollify the initial data,
\begin{equation*}
  v_{0\varepsilon}(\xi,s)=(J_{\varepsilon}*v_0)(\xi,s)
=\int_{\mathbb{R}}J_{\varepsilon}(\xi-y)v_0(y,s)dy
\in C([-\tau,0],W^{2,1}_w(\mathbb{R})\cap H^2(\mathbb{R})),
\end{equation*}
where $J_{\varepsilon}(\xi)$ is the mollifier. Let $v_{\varepsilon}(\xi,t)$ be the solution to (\ref{ch51}) with this mollified initial data. Using the similar arguments as the proof of Theorem \ref{thm-sta2}, we have
\begin{equation*}
  v_{\varepsilon}(\xi,t)\in C([0,\infty),W^{2,1}_w(\mathbb{R})\cap H^2(\mathbb{R})).
\end{equation*}
If $c>c^*$, we choose $\mu_1>0$ such that $\mathcal{M}(c,\mu_1)>0$, and if $c=c^*$, we can take only $\mu_1=0$ such that $\mathcal{M}(c^*,0)=\Delta(c^*,\lambda^*)=0$.  Multiplying (\ref{ch51}) by $w(\xi)e^{\mu_1s}$ and then integrating over $\mathbb{R}\times [0,t]$ with respect to $\xi$ and $t$, we have
\begin{equation}
  \label{ch53}
  \begin{split}
    &e^{\mu_1t}\int_{\mathbb{R}}w(\xi)v_{\varepsilon}(\xi,t)d\xi-d\int^t_0
    \int_{\mathbb{R}}e^{\mu_1s}w(\xi)\Delta_1v_{\varepsilon}(\xi,s)d\xi ds\\
    &\quad{}-\int^t_0\int_{\mathbb{R}}e^{\mu_1s}[\mu_1w(\xi)+cw'(\xi)
    +\partial_1f(0,0)w(\xi)]v_{\varepsilon}(\xi,s)d\xi ds\\
    &\quad{}-\partial_2f(0,0)\int^t_0\int_{\mathbb{R}}e^{\mu_1s}w(\xi)
    \left[\int_{\mathbb{R}}h(y)v_{\varepsilon}(\xi-y-c\tau,s-\tau)dy\right]d\xi ds
    \\
    &= \|v_0(0)\|_{L^1_{w}}+
    \int^t_0\int_{\mathbb{R}}e^{\mu_1s}w(\xi)Q(\xi,s)d\xi ds.
  \end{split}
\end{equation}
It is easy to see that, since $\lambda =-\frac{w'}{w}$,
\begin{multline}\label{ch54}
  \int^t_0\int_{\mathbb{R}}e^{\mu_1s}[\mu_1w(\xi)+cw'(\xi)
  +\partial_1f(0,0)w(\xi)]v_{\varepsilon}(\xi,s)d\xi ds\\
  =\int^t_0\int_{\mathbb{R}}e^{\mu_1s}[\mu_1-c\lambda+\partial_1f(0,0)]
  w(\xi)v_{\varepsilon}(\xi,s)d\xi ds
\end{multline}
and
\begin{equation}\label{ch55}
  \int^t_0\int_{\mathbb{R}}e^{\mu_1s}\Delta_1v_{\varepsilon}(\xi,s)w(\xi)d\xi ds=
  \int^t_0\int_{\mathbb{R}}e^{\mu_1s}[e^{\lambda}+e^{-\lambda}-2]w(\xi)v_{\varepsilon}(\xi,s)d\xi.
\end{equation}
By changing variables $\xi-y-c\tau\to \xi$, $s-\tau\to s$, and using the fact that
\begin{equation}\label{formu-h}
  \int_{\mathbb{R}}h(y)\frac{w(\xi+y+c\tau)}{w(\xi)}dy=e^{-\lambda
    c\tau}G(\lambda),
\end{equation}
we have
\begin{equation}\label{ch56}
\begin{split}
  &\int^t_0\int_{\mathbb{R}}e^{\mu_1s}w(\xi)
  \left[\int_{\mathbb{R}}h(y)v_{\varepsilon}(\xi-y-c\tau,s-\tau)dy\right]d\xi ds \\
  &=\int^{t-\tau}_{-\tau}\int_{\mathbb{R}}e^{\mu_1(s+\tau)}
  \left[\int_{\mathbb{R}}w(\xi+y+c\tau)h(y)dy\right]v_{\varepsilon}(\xi,s)d\xi ds \\
  &=e^{-\lambda
    c\tau}G(\lambda)\int^{t-\tau}_{-\tau}\int_{\mathbb{R}}e^{\mu_1(s+\tau)}
  w(\xi)v_{\varepsilon}(\xi,s)d\xi ds \\
  &\leq e^{-\lambda
    c\tau}G(\lambda)\int^{t}_{0}
  \int_{\mathbb{R}}e^{\mu_1(s+\tau)}w(\xi)v_{\varepsilon}(\xi,s)d\xi
  ds \\
  &\quad{}+e^{-\lambda
    c\tau}G(\lambda)\int^{0}_{-\tau}
  e^{\mu_1(s+\tau)}\|v_0(s)\|_{L^1_{w}(\mathbb{R})}ds.
\end{split}
\end{equation}
Applying  (\ref{ch54})--(\ref{ch56}) to (\ref{ch53}) and noticing that $Q(\xi,t) \leq 0$, we have
\begin{multline*}
  \label{ch58}   e^{\mu_1t}\int_{\mathbb{R}}w(\xi)v_{\varepsilon}(\xi,t)d\xi
  +\mathcal{M}(c,\mu_1)\int^t_0
  \int_{\mathbb{R}}e^{\mu_1s}w(\xi)v_{\varepsilon}(\xi,s)d\xi   ds\\ \leq
  \|v_0(0)\|_{L^1_{w}}
  +\partial_2f(0,0)e^{-\lambda
    c\tau}G(\lambda)\int^{0}_{-\tau}e^{\mu_1(s+\tau)}
  \|v_0(s)\|_{L^1_{w}(\mathbb{R})}ds.
\end{multline*}
It follows
that
\begin{equation*}
  \label{ch59}
  e^{\mu_1t}\int_{\mathbb{R}}w(\xi)v_{\varepsilon}(\xi,t)d\xi
  +\mathcal{M}(c,\mu_1)\int^t_0\int_{\mathbb{R}}e^{\mu_1s}
  w(\xi)v_{\varepsilon}(\xi,s)d\xi   ds\leq C.
\end{equation*}
We then establish the key energy estimate
\begin{equation}\label{ch591}
  e^{\mu_1t}\|v_{\varepsilon}(t)\|_{L^1_w(\mathbb{R})}
  +\int^t_0e^{\mu_1s}\|v_{\varepsilon}(s)\|_{L^1_w(\mathbb{R})}d\xi ds
  \leq C
\end{equation}
for $c>c^*$, where $\mathcal{M}(c,\mu_1)>0$, and
\begin{equation}\label{ch592}
  \|v_{\varepsilon}(t)\|_{L^1_w(\mathbb{R})}\leq C
\end{equation}
for $c=c^*$. Letting $\varepsilon\to 0$ in (\ref{ch591}) and (\ref{ch592}), we finally arrive at (\ref{ch52}) and (\ref{ch52'}). The proof is complete.  \qed

\begin{remark} \label{diff-Mei} Mei \emph{et al}. \cite{Mei4,Mei4a} assumed   that
  \begin{equation}\label{partial11}
    \mbox{$\partial_{11}f(u,v)<0$ for all $u$, $v\in [0,K]$},
  \end{equation}
  when they investigated the nonlinear stability of a continuous   reaction-diffusion equation.  In fact, under condition   (\ref{partial11}), we have $Q(\xi,t)\leq -C_1v^2(\xi,t)$ for some   positive constant $C_1$. Similarly, applying   (\ref{ch54})--(\ref{ch56}) to (\ref{ch53}) gives
  \begin{multline*}
    e^{\mu_1t}\int_{\mathbb{R}}w(\xi)v(\xi,t)d\xi
    +\mathcal{M}(c,\mu_1)\int^t_0\int_{\mathbb{R}}e^{\mu_1s}
    w(\xi)v(\xi,s)d\xi       ds\\
    +C_1\int^t_0e^{\mu_1s}\int_{\mathbb{R}}w(\xi)v^2(\xi,s)d\xi
    ds\leq C.
  \end{multline*}
  One then establishes the central energy estimate
  \begin{equation*}
    e^{\mu_1t}\|v(t)\|_{L^1_w(\mathbb{R})}
    +\int^t_0e^{\mu_1s}\|v(s)\|_{L^1_w(\mathbb{R})}d\xi ds
    +\int^{t}_0e^{\mu_1s}\|v(s)\|^2_{L^2_w(\mathbb{R})}ds\leq
    C
  \end{equation*}
  for $c>c^*$, and
  \begin{equation*}
    \|v(t)\|_{L^1_w(\mathbb{R})}+\int^{t}_0\|v(s)\|^2_{L^2_w(\mathbb{R})}ds\leq
    C
  \end{equation*}
  for $c=c^*$. This key energy estimate is much better than   (\ref{ch591}) and (\ref{ch592}) since it yields the boundedness of   $\int^{t}_0\|v(s)\|^2_{L^2_w(\mathbb{R})}ds$. Since we have a weaker   assumption, we need the following lemma to establish the boundedness   of $\int^{t}_0\|v(s)\|^2_{L^2_w(\mathbb{R})}ds$ without requiring   (\ref{partial11}).
\end{remark}

Although we cannot directly work with the original equations (\ref{ch51}) and (\ref{ch51'}) due to the lack of regularity of the solution, as discussed in the proof of Lemma \ref{lem-ch51}, we first consider a mollification, with solution
\begin{equation}\label{v-ep}
v_{\varepsilon}(\xi,t)\in C([0,\infty),H^{2}_w(\mathbb{R})\cap H^2(\mathbb{R})),
\end{equation}
and then take the limit $\varepsilon \to 0$ to obtain the corresponding energy estimate for the original solution $v(\xi,t)$. For the sake of simplicity, below we give the formal calculation using $v(\xi,t)$ directly to establish the desired energy estimates.

\begin{lemma}\label{lem-ch51-0} It holds that
  \begin{equation}\label{ch6-0}
    \|v(t)\|^2_{L^2_{w}(\mathbb{R})}
    +\int^{t}_0\|v(s)\|^2_{L^2_{w}(\mathbb{R})}ds\leq
    C
  \end{equation}
  for $c\geq c^*$, where $w(\xi)=e^{-\lambda(\xi-x_0)}$, $\lambda$ is   any fixed number in $(\lambda_1(c), \lambda^*)$ if $c>c^*$ but   $\lambda=\lambda^*$ if $c=c^*$.
\end{lemma}

\proof We multiply (\ref{ch51}) by $w(\xi)v(\xi,t)$ and notice that $Q(\xi,t)\leq 0$,
thus
\begin{equation}\label{ch6-1}
  (wv^2)_t+c(wv^2)_{\xi}+[c\lambda-2\partial_1f(0,0)]wv^2-d\cdot
  w\Delta_1v^2-2\partial_2f(0,0)wvh*v_{\tau}\leq 0.
\end{equation}
It follows from (\ref{v-ep}) that
$\lim_{\xi\to\pm\infty}w(\xi)v^2(\xi,t)=0$.
We then integrate (\ref{ch6-1}) over $\mathbb{R}\times [0,t]$ with respect to $\xi$ and $t$,
\begin{equation}\label{ch6-2}
\begin{split}
  &\|v(t)\|^2_{L^2_w(\mathbb{R})}-d\int^t_0
  \int_{\mathbb{R}}w(\xi)\Delta_1v^2(\xi,s)d\xi
  ds\\
  &+[c\lambda-2\partial_1f(0,0)]\int^t_0\int_{\mathbb{R}}w(\xi)v^2(\xi,s)d\xi
  ds\\
  &-2\partial_2f(0,0)\int^t_0\int_{\mathbb{R}}w(\xi)v(\xi,s)
  \left[\int_{\mathbb{R}}h(y)v(\xi-y-c\tau,s-\tau)dy\right]d\xi
  ds\\
  \leq &\|v_0(0)\|^2_{L^2_w(\mathbb{R})}.
  \end{split}
\end{equation}
Using Young's inequality $2|ab|\leq \eta a^2+\frac{1}{\eta}b^2$ for $\eta>0$ to be specified later, we obtain
\begin{equation}\label{ch6-3}
  \begin{split}
    &
    2\int^t_0\int_{\mathbb{R}}w(\xi)v(\xi,s)\left[\int_{\mathbb{R}}h(y)
      v(\xi-y-c\tau,s-\tau)dy\right]d\xi
    ds\\
    &=
    2\int^t_0\int_{\mathbb{R}}\int_{\mathbb{R}}h(y)w(\xi)v(\xi,s)
    v(\xi-y-c\tau,s-\tau)dyd\xi
    ds\\
    &\leq
    \eta\int^t_0\int_{\mathbb{R}}\int_{\mathbb{R}}h(y)w(\xi)v^2(\xi,s)dyd\xi\\
    &\quad{}+\frac{1}{\eta}\int^t_0\int_{\mathbb{R}}
    \int_{\mathbb{R}}h(y)w(\xi)v^2(\xi-y-c\tau,s-\tau)dyd\xi
    ds\\
    &=
    \eta\int^t_0\int_{\mathbb{R}}w(\xi)v^2(\xi,s)d\xi ds\\
    &\quad{}+\frac{1}{\eta}\int^{t-\tau}_{-\tau}\int_{\mathbb{R}}
    \int_{\mathbb{R}}h(y)w(\xi+y+c\tau)
    v^2(\xi,s)dyd\xi
    ds\\
    &\leq  \left[\eta+\frac{1}{\eta}e^{-\lambda
        c\tau}G(\lambda)\right]\int^t_0\|v(s)\|_{L^2_w(\mathbb{R})}ds
    +\frac{1}{\eta}e^{-\lambda
      c\tau}G(\lambda)\int^{0}_{-\tau}
    \|v(s)\|_{L^2_w(\mathbb{R})}
    ds.
  \end{split}
\end{equation}
Moreover,
\begin{equation}\label{ch6-4}
  \int_{\mathbb{R}}w(\xi)\triangle_1v^2(\xi,t)d\xi
  =(e^{\lambda}+e^{-\lambda}-2)\|v(t)\|^2_{L^2_w(\mathbb{R})}.
\end{equation}
An application of (\ref{ch6-3}) and (\ref{ch6-4}) to (\ref{ch6-2}) gives
\begin{equation}\label{ch6-5}
  \|v(t)\|^2_{L^2_w(\mathbb{R})}+A(\lambda,\eta)
  \int^t_0\|v(s)\|^2_{L^2_w(\mathbb{R})}ds\leq   C,
\end{equation}
where
\begin{multline*}
  A(\lambda,\eta)=\\
  \Delta(c,\lambda)-\partial_1f(0,0)+\partial_2f(0,0)e^{-\lambda     c\tau}G(\lambda)-\partial_2f(0,0)\left[\eta+\frac{1}{\eta}e^{-\lambda       c\tau}G(\lambda)\right].
\end{multline*}
In view of (F2) and (F3), we have $\partial_1f(0,0)\leq 0$ and $\partial_1f(0,0)+\partial_2f(0,0)>0$. Hence, $\partial_2f(0,0)e^{-\lambda c\tau}G(\lambda)-\partial_1f(0,0)>0$. Thus, $\Delta(c,\lambda)-\partial_1f(0,0)+\partial_2f(0,0)e^{-\lambda   c\tau}G(\lambda)>0$ for any fixed $\lambda\in (\lambda_1(c),\lambda^*)$ if $c>c^*$ but $\lambda=\lambda^*$ if $c=c^*$. Thus, it is easy to see that we can choose a suitable $\eta>0$ such that $A(\lambda,\eta)>0$. This proves (\ref{ch6-0}). \qed

Similarly to Lemma \ref{lem-ch51-0}, we have the following result.

\begin{lemma}\label{lem-ch51-01} It holds that
  \begin{equation}\label{ch6-01}
    \|v_{\xi}(t)\|^2_{L^2_{w}(\mathbb{R})}
    +\int^{t}_0\|v_{\xi}(s)\|^2_{L^2_{w}(\mathbb{R})}ds\leq
    C
  \end{equation}
  for $c\geq c^*$, where $w(\xi)=e^{-\lambda(\xi-x_0)}$, $\lambda$ is   any fixed number in $(\lambda_1(c), \lambda^*)$ if $c>c^*$ and   $\lambda=\lambda^*$ if $c=c^*$.
\end{lemma}

\proof Let us differentiate (\ref{ch51}) with respect to $\xi$, then we have
\begin{equation}\label{ch51-1}
    (v_{\xi})_t+cv_{\xi\xi}-d\cdot\Delta_1v_{\xi}
    -\partial_1f(0,0)v_{\xi}-\partial_2f(0,0)
    h*(v_{\tau})_{\xi}=Q_{\xi}.
\end{equation}
Thus, multiplying (\ref{ch51-1}) by $2w(\xi)v_{\xi}(\xi,t)$ yields
\begin{equation*}
\begin{split}
  &(wv_{\xi}^2)_t+c(wv_{\xi}^2)_{\xi}
  +[c\lambda-2\partial_1f(0,0)]wv_{\xi}^2-2d\cdot
  wv_{\xi}\Delta_1v_{\xi}\\
  &= 2\partial_2f(0,0)wv_{\xi}h*(v_{\tau})_{\xi}
 + 2wv_{\xi}Q_{\xi}.
    \end{split}
\end{equation*}
Notice that $2v_{\xi}\Delta_1v_{\xi}\leq \Delta_1v^2_{\xi}$ and
\begin{equation*}
  \begin{split}
    &
    v_{\xi}(\xi,t)Q_{\xi}(\xi,t)\\
    &=\left[\partial_1f\left(v+\phi,h*(v_{\tau}+\phi_{\tau})\right)
    -\partial_1f\left(\phi,h*(\phi_{\tau})\right) \right]\phi_{\xi}v_{\xi}\\
    &\quad{}+\left[\partial_2f\left(v+\phi,h*(v_{\tau}+\phi_{\tau})\right)
    -\partial_2f\left(\phi,h*(\phi_{\tau})\right) \right]h*(\phi_{\tau})_{\xi}v_{\xi}\\
    &\quad{}+\left[\partial_1f\left(v+\phi,h*(v_{\tau}+\phi_{\tau})\right)
    -\partial_1f\left(0,0)\right) \right]v^2_{\xi}\\
    &\quad{}+\left[\partial_2f\left(v+\phi,h*(v_{\tau}+\phi_{\tau})\right)
    -\partial_2f\left(0,0)\right) \right]h*(v_{\tau})_{\xi}v_{\xi}\\
    &\leq C_1\left[\phi_{\xi}+h*(|(v_{\tau})_{\xi}|
    +(\phi_{\tau})_{\xi})\right]|v_{\xi}|.
  \end{split}
\end{equation*}
Then, using Young's inequality $2|ab|\leq \eta a^2+\frac{1}{\eta}b^2$ for $\eta>0$ to be specified later, we obtain
\begin{equation}\label{ch6-1-1}
\begin{split}
  &(wv_{\xi}^2)_t+c(wv_{\xi}^2)_{\xi}
  +[c\lambda-2\partial_1f(0,0)]wv_{\xi}^2-d\cdot
  w\Delta_1v_{\xi}^2\\
  &\leq 2\partial_2f(0,0)wv_{\xi}h*(v_{\tau})_{\xi}
 + 2C_1w(\xi)\left[\phi_{\xi}+h*(|(v_{\tau})_{\xi}|
    +(\phi_{\tau})_{\xi})\right]|v_{\xi}|\\
    &\leq 2C_1w(\xi)\left[\phi_{\xi}
    +h*(\phi_{\tau})_{\xi}\right]|v_{\xi}|+2C_2w(\xi)|v_{\xi}|h*(|(v_{\tau})_{\xi}|)\\
    &\leq w(\xi)\left[\frac{C_1}{\eta}\phi^2_{\xi}
    +(2C_1+C_2)\eta v^2_{\xi}+\frac{C_2}{\eta}h*(|(v_{\tau})_{\xi}|)^2+
    \frac{C_1}{\eta}h*((\phi_{\tau})_{\xi})^2\right].
    \end{split}
\end{equation}
where $C_2=C_1+\partial_2f(0,0)>0$.
We then integrate (\ref{ch6-1-1}) over $\mathbb{R}\times [0,t]$ with respect to $\xi$ and $t$,
\begin{equation}\label{ch6-2-1}
\begin{split}
  &\|v_{\xi}(t)\|^2_{L^2_w(\mathbb{R})}-d\int^t_0
  \int_{\mathbb{R}}w(\xi)\Delta_1v_{\xi}^2(\xi,s)d\xi
  ds\\
  &\quad{}
  +[c\lambda-2\partial_1f(0,0)]
  \int^t_0\int_{\mathbb{R}}w(\xi)v_{\xi}^2(\xi,s)d\xi
  ds \\
  &\leq \|v_0(0)\|^2_{L^2_w(\mathbb{R})}+
  \left[(2C_1+C_2)\eta
  +\frac{C_2}{\eta}e^{-c\tau\lambda}G(\lambda)\right] \int^t_{0}\|v_{\xi}(s)\|^2_{L^2_w(\mathbb{R})}ds\\
  &\quad{}+\frac{C_1}{\eta}\left[1+e^{-c\tau\lambda}G(\lambda)\right]\int^t_{0}\|\phi_{\xi}\|^2_{L^2_w(\mathbb{R})}ds\\
  &\quad{}+\frac{C_2}{\eta}
  e^{-c\tau\lambda}G(\lambda)\int^0_{-\tau}\|v_{\xi}(s)\|^2_{L^2_w(\mathbb{R})}
  ds.
\end{split}
\end{equation}
In view of Theorems \ref{thm-as} and \ref{thm-as2}, we have $\|\phi_{\xi}\|^2_{L_w^2(\mathbb{R})}<C$, and hence
\begin{equation}\label{ch6-5-1}
  \|v_{\xi}(t)\|^2_{L^2_w(\mathbb{R})}+B(\lambda,\eta)
  \int^t_0\|v_{\xi}(s)\|^2_{L^2_w(\mathbb{R})}ds\leq   C,
\end{equation}
where
\begin{multline*}
  B(\lambda,\eta)=\\
  \Delta(c,\lambda)-\partial_1f(0,0)+\partial_2f(0,0)e^{-\lambda     c\tau}G(\lambda)-\left[(2C_1+C_2)\eta
  +\frac{C_2}{\eta}e^{-c\tau\lambda}G(\lambda)\right].
\end{multline*}
Using similar arguments as those of the proof of Lemma \ref{lem-ch51-0}, it is easy to see that we can choose a suitable $\eta>0$ such that $B(\lambda,\eta)>0$. This proves (\ref{ch6-01}). \qed

\begin{lemma}\label{lem27}
  For any $c\geq c^*$, it holds that
  \begin{equation}\label{ch510}
    \|v(t)\|^2_{L^2(\mathbb{R})}+\int^{t}_0\|v(s)\|^2_{L^2(\mathbb{R})}ds\leq
    C.
  \end{equation}
\end{lemma}

\proof In view of $\partial_1f(K,K)+\partial_2f(K,K)<0$ (by (F2)) and $\partial_2f(K,K)\geq 0$ (by (F1)), we have $\partial_1f(K,K)<-\partial_2f(K,K)\leq0$.  Hence, the equation $\partial_2f(K,K)x^2+2\partial_1(K,K)x+\partial_2f(K,K)=0$ has two distinct positive solutions
\begin{equation*}
  x_{\pm}=\frac{-\partial_1f(K,K)\pm
    \sqrt{[\partial_1f(K,K)]^2-[\partial_2f(K,K)]^2}}{\partial_2f(K,K)}.
\end{equation*}
Moreover, $\partial_2f(K,K)x^2+2\partial_1(K,K)x+\partial_2f(K,K)<0$ for all $x\in (x_-,x_+)$. It follows that we can choose a positive constant $\eta\in (x_-,x_+)$ such that
\begin{equation}\label{ch516}
  2\partial_1f(K,K)+\eta\partial_2f(K,K)+\frac{1}{\eta}\partial_2f(K,K)
  <0.
\end{equation}
Thus, (\ref{ch517*}) and (\ref{ch516}) implies that we can choose $x_0$ sufficiently large such that
\begin{equation}\label{ch518}
  C_3=-2G_1(x_0)-\eta G_2(x_0)-\frac{1}{\eta}B(x_0)>0.
\end{equation}
In view of $G_j(\xi)\leq G_j(x_0)$ for $\xi\geq x_0$, $j=1,2$, we have
\begin{equation}\label{ch519}
  -2G_1(\xi)-\eta G_2(\xi)-\frac{1}{\eta}B(x_0)\geq C_3>0.
\end{equation}
Moreover, since $w(\xi)=e^{-\lambda(x-x_0)}\geq 1$ for all $x\in (-\infty,x_0]$, Lemma~\ref{lem-ch51-0} guarantees that
\begin{equation}\label{ch511}
\int^{\infty}_0\int^{x_0}_{-\infty}v^2(\xi,s)d\xi ds<C.
\end{equation}

We rewrite (\ref{ch51}) as
\begin{equation}\label{ch51*2}
  v_t+cv_{\xi}-d\cdot \Delta_1v-G_1v-G_2h*v_{\tau}=\widetilde{Q},
\end{equation}
where $v=v(\xi,t)$, $v_{\tau}=v(\xi-c\tau,t-\tau)$, $G_1=G_1(\xi)$, $G_2=G_2(\xi)$, $\widetilde{Q}=f(v+\phi,h*(v_{\tau}+\phi_{\tau}))-f(\phi,h*\phi_{\tau})-G_1(\xi)v-G_2(\xi)h*v_{\tau}$, $\phi=\phi(\xi)$, $\phi_{\tau}=\phi(\xi-c\tau)$.  Clearly, $\widetilde{Q}(\xi,t)\leq 0$.  A multiplication of (\ref{ch51*2}) by $v(\xi,t)$ gives then
\begin{equation*}
  (v^2)_t+c(v^2)_{\xi}-2dv\Delta_1v-2G_1v^2-2G_2vh*v_{\tau}\leq 0.
\end{equation*}
Integrating this over $\mathbb{R}\times [0,t]$ with respect to $\xi$ and $t$ and applying inequality (\ref{ch6-6}), we obtain
\begin{multline}\label{ch512}
  \|v(t)\|^2_{L^2(\mathbb{R})}-2\int^t_0
  \int_{\mathbb{R}}G_1(\xi)v^2(\xi,s)d\xi ds\\
  -2\int^t_0\int_{\mathbb{R}}G_2(\xi)v(\xi,s)
  \left[\int_{\mathbb{R}}h(y)v(\xi-y-c\tau,s-\tau)dy\right]d\xi
  ds\leq \|v_0(0)\|^2_{L^2(\mathbb{R})}.
\end{multline}
From Young's inequality $2|ab|\leq \eta a^2+\frac{1}{\eta}b^2$ for $\eta\in (x_-,x_+)$, we obtain
\begin{equation}\label{ch513}
  \begin{split}
    &
    2\int^t_0\int_{\mathbb{R}}G_2(\xi)v(\xi,s)
    \left[\int_{\mathbb{R}}h(y)v(\xi-y-c\tau,s-\tau)dy\right]d\xi
    ds\\
    &=
    2\int^t_0\int_{\mathbb{R}}\int_{\mathbb{R}}h(y)
    G_2(\xi)v(\xi,s)v(\xi-y-c\tau,s-\tau)dyd\xi ds\\
    &\leq
    \eta\int^t_0\int_{\mathbb{R}}\int_{\mathbb{R}}h(y)G_2(\xi)v^2(\xi,s)dy
    d\xi\\
    &\quad{}+\frac{1}{\eta}\int^t_0\int_{\mathbb{R}}
    \int_{\mathbb{R}}h(y)G_2(\xi)v^2(\xi-y-c\tau,s-\tau)dyd\xi
    ds\\
    &=
    \eta\int^t_0\int_{\mathbb{R}}G_2(\xi)v^2(\xi,s)d\xi ds\\
    &\quad{}+\frac{1}{\eta}\int^{t-\tau}_{-\tau}\int_{\mathbb{R}}
    \int_{\mathbb{R}}h(y)G_2(\xi+y+c\tau)
    v^2(\xi,s)dyd\xi
    ds\\
    &\leq \eta\partial_2f\left(0,0\right)\int^t_0
    \int^{x_0}_{-\infty}v^2(\xi,s)d\xi ds
    +\eta\int^t_0\int^{\infty}_{x_0}G_2(\xi)v^2(\xi,s)d\xi ds\\
    &\quad{}+\frac{1}{\eta}\int^{t-\tau}_{-\tau}
    \int_{\mathbb{R}}B(\xi)v^2(\xi,s)d\xi ds\\
    &\leq
    C+\eta\int^t_0\int^{\infty}_{x_0}G_2(\xi)v^2(\xi,s)d\xi ds
    +\frac{1}{\eta}B(\xi_0)\int^{t}_{0}\int^{\infty}_{x_0}v^2(\xi,s)d\xi
    ds.
  \end{split}
\end{equation}
An application of (\ref{ch513}) to (\ref{ch512}) gives
\begin{multline}\label{ch515}
  \|v(t)\|^2_{L^2(\mathbb{R})}-2\int^t_0\int^{x_0}_{-\infty}
  G_1(\xi)v^2(\xi,s)d\xi ds\\
  -\int^t_0\int^{\infty}_{x_0}\left[2G_1(\xi)
    +\eta G_2(\xi)+\frac{1}{\eta}B(\xi_0)\right]v^2(\xi,s)d\xi
  ds\leq   C.
\end{multline}
Then (\ref{ch519}) yields for (\ref{ch515}), with $C_3\int^t_0\int^{x_0}_{-\infty}v^2(\xi,s)d\xi ds$ added,
\begin{equation}\label{ch520}
  \|v(t)\|^2_{L^2(\mathbb{R})}+C_3\int^t_0\|v(s)\|^2_{L^2(\mathbb{R})}ds
  \leq     C,
\end{equation}
which proves (\ref{ch510}).  \qed

\begin{remark}
For the partial differential equation (\ref{eqcon1}), we have
\begin{equation}\label{ch51*2-11}
  v_t+cv_{\xi}-dv_{\xi\xi}-G_1v-G_2h*v_{\tau}=\widetilde{Q},
\end{equation}
with the initial data $v(\xi,s)=v_0(\xi,s)>0$, $s\in [-\tau,0]$, where $v$, $v_{\tau}$, $\phi_{\tau}$, and $\widetilde{Q}$ are defined as for (\ref{ch51*2}).  We can derive the following $L^2$-energy estimate for $v(\xi, t)$ (see \cite{LvWang,Mei4,Mei4a} for more details)
\begin{equation}\label{ch510b}
    \|v(t)\|^2_{L^2(\mathbb{R})}+\int^{t}_0\|v(s)\|^2_{H^1(\mathbb{R})}ds\leq
    C, \quad t\geq 0.
\end{equation}
It is important to notice the difference in the integrands in (\ref{ch510}) and (\ref{ch510b}) (weaker norm for the difference-differential equation considered here). In order to derive for the partial differential equation (\ref{eqcon1}) the $L^2$-energy estimate for $v_{\xi}(\xi,t)$, we can differentiate (\ref{ch51*2-11}) with respect to $\xi$, multiply the resulting equation by $v_{\xi}(\xi,t)$ and then integrate it over $\mathbb{R}\times [0,t]$ with respect to $\xi$ and $t$. Employing the key estimates (\ref{ch510b}), we obtain an $L^2$-energy estimate for $v_{\xi}(\xi,t)$ (see Lemma 3.3 in \cite{Mei4}):
\begin{equation*}
    \|v_{\xi}(t)\|^2_{L^2(\mathbb{R})}
    +\int^{t}_0\|v_{\xi}(s)\|^2_{H^1(\mathbb{R})}ds\leq
    C, \quad t\geq 0.
\end{equation*}
In contrast, for equation (\ref{ch51*2}), the $L^2$-energy estimate (\ref{ch510}) cannot help us to obtain the $L^2$-energy estimate for $v_{\xi}(\xi,t)$.
\end{remark}

Next, we derive an $L^2$-energy estimate for $v_{\xi}(\xi,t)$.
\begin{lemma}
  For any $c\geq c^*$, it holds that
  \begin{equation}\label{ch524}
    \|v_{\xi}(t)\|^2_{L^2(\mathbb{R})}
    +\int^{t}_0\|v_{\xi}(s)\|^2_{L^2(\mathbb{R})}ds\leq
    C.
\end{equation}
\end{lemma}

\proof Consider the function
\begin{equation*}
  H(\xi)=\partial_1f(K,K)+3\partial_2f(K,K)
-5G_1(\xi)-4G_2(\xi)-3B(\xi)
\end{equation*}
for $\xi\in \mathbb{R}$. It is easy to see that $H$ is continuous and bounded on $\mathbb{R}$.  Note that $H(+\infty)=-4[\partial_1f(K,K)+\partial_2f(K,K)]>0$. Thus there exists $x_0\in \mathbb{R}$ such that
\begin{equation*}
  H(\xi)\geq H(x_0)>0
\end{equation*}
for all $\xi\geq x_0$. Moreover, since $w(\xi)=e^{-\lambda(x-x_0)}\geq 1$ for all $x\in (-\infty,x_0]$, Lemma~\ref{lem-ch51-01} guarantees that
\begin{equation}\label{ch511-1}
\int^{\infty}_0\int^{x_0}_{-\infty}v_{\xi}^2(\xi,s)d\xi ds<C.
\end{equation}

Let us differentiate (\ref{ch51*2}) with
respect to $\xi$ and multiply the resulting equation by $v_{\xi}(\xi,t)$. We then obtain
\begin{equation}\label{ch521}
  \begin{split}
    &
    (v^2_{\xi})_t+c(v^2_{\xi})_{\xi}-2dv_{\xi}[v_{\xi}(\xi+1,t)+v_{\xi}
    (\xi-1,t)-2v_{\xi}(\xi,t)]\\
    &\quad{}-2G_1(\xi)v^2_{\xi}-2G_2(\xi)v_{\xi}\int_{\mathbb{R}}h(y)v_{\xi}
    (\xi-y-c\tau,s-\tau)dy\\
    &=
    2v_{\xi}\widetilde{Q}_{\xi}+2G'_1(\xi)vv_{\xi}+2G'_2(\xi)
    v_{\xi}\int_{\mathbb{R}}h(y)v(\xi-y-c\tau,s-\tau)dy.
  \end{split}
\end{equation}
Integrating (\ref{ch521}) over $\mathbb{R}\times [0, t]$ with respect to $\xi$ and $t$, we obtain
\begin{equation}\label{ch522}
  \begin{split}
    &   \|v_{\xi}(t)\|^2_{L^2(\mathbb{R})}
    -2\int^t_0\int_{\mathbb{R}}G_1(\xi)v^2_{\xi}(\xi,s)d\xi ds\\
    &\quad{}-2\int^t_0\int_{\mathbb{R}}G_2(\xi)v_{\xi}(\xi,s)\int_{\mathbb{R}}h(y)v_{\xi}
    (\xi-y-c\tau,s-\tau)dyd\xi ds\\
    &\leq  C+
    2\int^t_0\int_{\mathbb{R}}\left[v_{\xi}(\xi,s)\widetilde{Q}_{\xi}(\xi,s)
      +G'_1(\xi)v(\xi,s)v_{\xi}(\xi,s)\right.\\
    &\left.\quad{}+G'_2(\xi)v_{\xi}(\xi,s)\int_{\mathbb{R}}h(y)
      v(\xi-y-c\tau,s-\tau)dy
    \right]d\xi ds.
  \end{split}
\end{equation}
Now let us examine the inner integral of the last term on the right-hand side of (\ref{ch522}).  In view of the expression for $\widetilde{Q}(\xi,t)$, we have
\begin{equation*}
  \begin{split}
    &\widetilde{Q}_{\xi}+G'_1(\xi)v+G'_2(\xi)
    h*v_{\tau}\\
    &=
    \left\{f(v+\phi,h*(v_{\tau}+\phi_{\tau})))-f(\phi,h*\phi_{\tau})
    \right\}_{\xi}
    -G_1(\xi)v_{\xi}-G_2(\xi)h*(v_{\tau})_{\xi}\\
    &=  \left[\partial_1f(v+\phi,h*(v_{\tau}+\phi_{\tau}))
      -\partial_1f(\phi,h*\phi_{\tau})\right][v_{\xi}+\phi'(\xi)]\\
    &\quad{}+\left[\partial_2f(v+\phi,h*(v_{\tau}+\phi_{\tau}))
      -\partial_2f(\phi,h*\phi_{\tau})\right][h*(v_{_{\tau}}+\phi_{_{\tau}})_{\xi}].
  \end{split}
\end{equation*}
In view of (\ref{ch517}) and assumption (F3), we have
\begin{equation*}
  \partial_1f(K,K)
  -G_1(\xi)\leq \partial_1f(v+\phi,h*(v_{\tau}+\phi_{\tau}))
  -\partial_1f(\phi,h*\phi_{\tau})\leq 0
\end{equation*}
and
\begin{equation*}
  \partial_2f(K,K)
  -G_2(\xi)\leq   \partial_2f(v+\phi,h*(v_{\tau}+\phi_{\tau}))
  -\partial_2f(\phi,h*\phi_{\tau})\leq 0.
\end{equation*}
Thus, using Young's inequality, we have
\begin{equation*}\label{ch525}
\begin{split}
  & \int_{\mathbb{R}}\left[\widetilde{Q}_{\xi}+G'_1(\xi)v+G'_2(\xi)
    h*v_{\tau}\right]v_{\xi}d\xi\\
  &\leq
  \int_{\mathbb{R}}[G_1(\xi)-\partial_1f(K,K)]|\phi'(\xi)|
    \left|v_{\xi}\right|d\xi\\
  &\quad{}+\int_{\mathbb{R}}[G_2(\xi)-\partial_2f(K,K)]
  |h*(v_{\tau}+\phi_{\tau})_{\xi}|
  |v_{\xi}|d\xi\\
  &\leq
  \frac{1}{2}\int_{\mathbb{R}}[G_1(\xi)-\partial_1f(K,K)]|\phi'(\xi)|^2d\xi
  +
    \frac{1}{2}\int_{\mathbb{R}}[G_1(\xi)-\partial_1f(K,K)]
    \left|v_{\xi}(\xi,t)\right|^2d\xi\\
  &\quad{}+\int_{\mathbb{R}}[G_2(\xi)-\partial_2f(K,K)]
  |h*(v_{\tau})_{\xi}||v_{\xi}(\xi,t)|d\xi\\
  &\quad{}+\frac{1}{2}
  \int_{\mathbb{R}}[G_2(\xi)-\partial_2f(K,K)]
  (\phi_{\tau})^2_{\xi}|d\xi\\
  &\quad{}+\frac{1}{2}
  \int_{\mathbb{R}}[G_2(\xi)-\partial_2f(K,K)]
  |v_{\xi}(\xi,t)|^2d\xi\\
  &\leq
  C+\frac{1}{2}\int_{\mathbb{R}}[G_1(\xi)+G_2(\xi)
  -\partial_1f(K,K)-\partial_2f(K,K)]
    \left|v_{\xi}(\xi,t)\right|^2d\xi\\
  &\quad{}+\int_{\mathbb{R}}[G_2(\xi)-\partial_2f(K,K)]
  |h*(v_{\tau})_{\xi}||v_{\xi}(\xi,t)|d\xi.
\end{split}
\end{equation*}
This, together with (\ref{ch522}), implies that
\begin{equation*}
  \begin{split}
    &   \|v_{\xi}(t)\|^2_{L^2(\mathbb{R})}
    -\frac{1}{2}\int^t_0\int_{\mathbb{R}}\left[5G_1(\xi)+G_2(\xi)
    -\partial_1f(K,K)-\partial_2f(K,K)\right]v^2_{\xi}(\xi,s)d\xi ds\\
    &\leq  C+\int^t_0\int_{\mathbb{R}}[3G_2(\xi)-\partial_2f(K,K)]|v_{\xi}(\xi,s)|\int_{\mathbb{R}}h(y)|v_{\xi}
    (\xi-y-c\tau,s-\tau)|dyd\xi ds.
\end{split}
\end{equation*}
Using a similar argument as that in the proof of Lemma \ref{lem27}, we have
$$
\|v_{\xi}(t)\|^2_{L^2(\mathbb{R})}+
    \int^t_0\int_{\mathbb{R}}H(\xi)v^2_{\xi}(\xi,s)d\xi ds\leq   C.
$$
Namely,
$$
\|v_{\xi}(t)\|^2_{L^2(\mathbb{R})}+
    \int^t_0\int^{x_0}_{-\infty}H(\xi)v^2_{\xi}(\xi,s)d\xi ds+
    H(x_0)\int^t_0\int_{x_0}^{+\infty}v^2_{\xi}(\xi,s)d\xi ds
    \leq   C.
$$
Recall that $H(\cdot)$ is continuous and bounded on $\mathbb{R}$, then we have
$$
\|v_{\xi}(t)\|^2_{L^2(\mathbb{R})}+
    H(x_0)\int^t_0\int_{x_0}^{+\infty}v^2_{\xi}(\xi,s)d\xi ds
    \leq   C.
$$
This, together with (\ref{ch511-1}), implies that we obtain the first order estimate (\ref{ch524}).  \qed

Based on the above lemmas, we can prove the following two convergence results.  One is the exponential stability for the noncritical traveling waves with $c > c^*$, and the other one is the algebraic stability for the critical traveling wave with $c = c^*$. We first prove the exponential stability.
\begin{lemma}
  For any $c> c^*$, it holds that
  \begin{equation}\label{ch525-2}
    \sup\limits_{x\in \mathbb{R}}|u^+(x,t)-\phi(x+ct)|
    =\|v(t)\|_{L^{\infty}(\mathbb{R})}\leq
    Ce^{-\mu_2t/3},\quad t\geq 0,
  \end{equation}
  where $\mu_2\in (0, \mu_1]$ satisfies $\mathcal{N}(\mu_2)<0$.
\end{lemma}

\proof For the choice of $x_0$ in (\ref{assume2}), it is easy to see that
\begin{equation*}
  \|v(t)\|^2_{L^{2}(-\infty,x_0]}=\int^{x_0}_{-\infty}|v(\xi,t)|^2d\xi\leq
  \|v(t)\|_{L^{\infty}(-\infty,x_0]}\|v(t)\|_{L^{1}(-\infty,x_0]}
\end{equation*}
and
\begin{multline*}
  v^2(\xi,t)=\int^{\xi}_{-\infty}[v^2(\xi,t)]_{\xi}d\xi
  =2\int^{\xi}_{-\infty}v(\xi,t)v_{\xi}(\xi,t)d\xi \\
  \leq
  2\|v(t)\|_{L^{2}(-\infty,x_0]}\|v_{\xi}(t)\|_{L^{2}(-\infty,x_0]}
\end{multline*}
for all $\xi\in (-\infty,x_0]$.  Thus we have
\begin{equation*}
  \|v(t)\|^2_{L^{\infty}(-\infty,x_0]}\leq
  2\|v(t)\|^{\frac{1}{2}}_{L^{\infty}(-\infty,x_0]}
  \|v(t)\|^{\frac{1}{2}}_{L^{1}(-\infty,x_0]}\|v_{\xi}(t)\|_{L^{2}(-\infty,x_0]},
\end{equation*}
and so
\begin{equation}\label{ch526}
  \|v(t)\|_{L^{\infty}(-\infty,x_0]}\leq
  \sqrt[3]{4}\|v(t)\|^{\frac{1}{3}}_{L^{1}(-\infty,x_0]}
  \|v_{\xi}(t)\|^{\frac{2}{3}}_{L^{2}(-\infty,x_0]}.
\end{equation}
Since $w(\xi)=e^{-\lambda(x-x_0)}\geq 1$ for all $\xi\in (-\infty, x_0]$, it follows from (\ref{ch52}) that
\begin{equation*}
  \|v(t)\|_{L^{1}(-\infty,x_0]}\leq   \|v(t)\|_{L^{1}_{w}(-\infty,x_0]}\leq Ce^{-\mu_1t}
\end{equation*}
for all $t\geq 0$. This, together with (\ref{ch524}) and (\ref{ch526}), implies that
\begin{equation}\label{ch527-2}
  \|v(t)\|_{L^{\infty}(-\infty,x_0]}\leq
  Ce^{-\mu_1t/3},\quad t\geq 0.
\end{equation}

A multiplication of (\ref{ch51}) by $e^{\mu_2t}$ and an integration over $\mathbb{R}\times [0,t]$ with respect to $\xi$ and $t$ combined with the fact that $\widetilde{Q}(\xi,t)\leq 0$ yield
\begin{multline}\label{ch528-2}
  e^{\mu_2t}\|v(t)\|_{L^1(\mathbb{R})}
  -\int^t_0\int_{\mathbb{R}}e^{\mu_2s}[\mu_2+G_1(\xi)]v(\xi,s)d\xi ds
  -d\int^t_0\int_{\mathbb{R}}e^{\mu_2s}\Delta_1v(\xi,s)d\xi ds  \\
  \leq    \|v_0(0)\|_{L^1(\mathbb{R})}
  +\int^t_0\int_{\mathbb{R}}e^{\mu_2s}G_2(\xi)
  \left[\int_{\mathbb{R}}h(y)v(\xi-y-c\tau,s-\tau)dy\right]d\xi ds.
\end{multline}
By the change of variables $\xi-y-c\tau\to \xi$ and $s-\tau\to s$, we have
\begin{equation}\label{ch529}
  \begin{split}
    &\int^t_0\int_{\mathbb{R}}e^{\mu_2s}G_2(\xi)
    \left[\int_{\mathbb{R}}h(y)v(\xi-y-c\tau,s-\tau)dy\right]d\xi ds
    \\
    &=     \int^{t-\tau}_{-\tau}\int_{\mathbb{R}}e^{\mu_2(s+\tau)}
    \left[\int_{\mathbb{R}}h(y)G_2(\xi+y+c\tau)dy\right]
    v(\xi,s)d\xi     ds \\
    &=     e^{\mu_2\tau}\int^{t-\tau}_{-\tau}
    \int_{\mathbb{R}}e^{\mu_2s}B(\xi)v(\xi,s)d\xi ds.
  \end{split}
\end{equation}
Substituting (\ref{ch529}) into (\ref{ch528-2}) and noticing that $\int_{\mathbb{R}}e^{\mu_2s}\Delta_1v(\xi,s)d\xi=0$, we obtain
\begin{multline}\label{ch530}
  e^{\mu_2t}\|v(t)\|_{L^1(\mathbb{R})}-\int^t_0
  \int_{\mathbb{R}}e^{\mu_2s}[\mu_2+G_1(\xi)]v(\xi,s)d\xi ds \\
  \leq
  \|v_0(0)\|_{L^1(\mathbb{R})}+e^{\mu_2\tau}\int^{t-\tau}_{-\tau}
  \int_{\mathbb{R}}e^{\mu_2s}B(\xi)v(\xi,s)d\xi ds.
\end{multline}
Splitting each integral of (\ref{ch530}) into two parts $(-\infty,x_0]\cup [x_0,\infty)$, we have
\begin{multline}\label{ch531}
  e^{\mu_2t}\|v(t)\|_{L^1[x_0,\infty)}
  -\int^t_0\int_{x_0}^{\infty}e^{\mu_2s}[\mu_2+G_1(\xi)
  +e^{\mu_2\tau}B(\xi)]v(\xi,s)d\xi ds \\
  \leq    \|v_0(0)\|_{L^1(\mathbb{R})}-J(t)+e^{\mu_2\tau}
  \int^{0}_{-\tau}\int_{x_0}^{\infty}e^{\mu_2s}B(\xi)v(\xi,s)d\xi ds,
\end{multline}
where
\begin{multline*}
  J(t)=e^{\mu_2t}\|v(t)\|_{L^1(-\infty,x_0]}
  -\int^t_0\int^{x_0}_{-\infty}e^{\mu_2s}[\mu_2+G_1(\xi)]
  v(\xi,s)d\xi ds\\
  -e^{\mu_2\tau}\int^{t-\tau}_{-\tau}\int^{x_0}_{-\infty}
  e^{\mu_2s}B(\xi)v(\xi,s)d\xi ds.
\end{multline*}
Since $G_1(\xi)$ and $B(\xi)$ are both non-increasing, $\partial_1f(0,0)\geq G_1(\xi)\geq G_1(x_0)$ and $\partial_2f(0,0)\geq B(\xi)\geq B(x_0)$ for all $\xi\in (-\infty,x_0]$.  Moreover, in view of Lemma \ref{lem-ch51} and $w(\xi)\geq 1$ for all $\xi\in (-\infty,x_0]$, it follows that
\begin{equation*}
  e^{\mu_2t}\|v(t)\|_{L^1(-\infty,x_0]}<C,\quad
  \int^t_0e^{\mu_2s}\|v(s)\|_{L^1(-\infty,x_0]} ds<C.
\end{equation*}
This means that $J(t)$ is bounded. Thus, it follows from (\ref{assume2}) and (\ref{ch531}) that
\begin{equation*}
  e^{\mu_2t}\|v(t)\|_{L^1[x_0,\infty)}\leq  C.
\end{equation*}
By a similar argument as above, we obtain
\begin{equation}\label{ch532}
  \|v(t)\|_{L^{\infty}[x_0,\infty)}
  \leq\sqrt[3]{4}\|v(t)\|^{\frac{1}{3}}_{L^{1}[x_0,\infty)}
  \|v_{\xi}(t)\|^{\frac{2}{3}}_{L^{2}[x_0,\infty)}\leq Ce^{-\mu_2t/3}.
\end{equation}
Combining (\ref{ch527-2}) and (\ref{ch532}) gives (\ref{ch525-2}). This completes the proof.  \qed

Now we are going to prove the algebraic stability for the critical traveling wave with $c = c^*$. In this case, we have $\lambda=\lambda_1(c)=\lambda^*$. Using the linearization of (\ref{ch51}) at 0 and noticing that $Q(\xi,t)\leq 0$, we find
\begin{multline}\label{ch51*}
  v_t+c^*v_{\xi}-d[v(\xi+1,t)+v(\xi-1,t)-2v(\xi,t)] \\
  \partial_1f(0,0)v-\partial_2f(0,0)\int_{\mathbb{R}}
  h(y)v(\xi-y-c^*\tau,t-\tau)dy   \leq 0.
\end{multline}
Let $\bar{v}(\xi,t)$ be the solution of the equation
\begin{multline}\label{ch52*}
  v_t+c^*v_{\xi}-d[v(\xi+1,t)+v(\xi-1,t)-2v(\xi,t)] \\
  -\partial_1f(0,0)v-\partial_2f(0,0)\int_{\mathbb{R}}
  h(y)v(\xi-y-c^*\tau,t-\tau)dy= 0
\end{multline}
with the initial data $v_0(\xi,s)$, $s\in [-\tau,0]$. By the comparison principle
\begin{equation*}
  v(\xi,t)\leq \bar{v}(\xi,t)\quad\mbox{for all $(\xi,t)\in
    \mathbb{R}\times \mathbb{R}^+$}.
\end{equation*}
Let $\tilde{v}(\xi,t)=w(\xi)\bar{v}(\xi,t)$, where $w(\xi)=\exp\{-\lambda^*(\xi-x_0)\}$ and $\lambda^*$ is defined in Lemma \ref{lem2}. Then
\begin{multline}\label{tiv}
  \tilde{v}_t+c^*\tilde{v}_{\xi}-d[e^{\lambda^*}\tilde{v}(\xi+1,t)
  +e^{-\lambda^*}\tilde{v}(\xi-1,t)-2\tilde{v}(\xi,t)]
  +[c^*\lambda^*-\partial_1f(0,0)]\tilde{v} \\
  -\partial_2f(0,0)\int_{\mathbb{R}}h(y)\exp\{-\lambda^*(y+c^*\tau)\}
  \tilde{v}(\xi-y-c^*\tau,t-\tau)dy=0.
\end{multline}

\begin{lemma}\label{lem-tiv}
  $\|\tilde{v}(t)\|_{L^{\infty}(\mathbb{R})}\leq Ct^{-1/2}$ for all   $t>0$.
\end{lemma}

\proof We take the Fourier transform of equation (\ref{tiv}) with respect to $\xi$ and write the Fourier transforms of $\hat{v}(\xi,t)$ as $V(\omega,t)$. Then
\begin{equation}\label{eqV}
  V_t(\omega,t)+c_1(\omega)V(\omega,t)=c_2(\omega)V(\omega,t-\tau),
\end{equation}
where
\begin{align*}
  c_1(\omega)&=   ic^*\omega+c^*\lambda^*-\partial_1f(0,0)
  -d[\exp\{\lambda^*+i\omega\}+\exp\{-\lambda^*-i\omega\}-2], \\
  c_2(\omega)&=
  \partial_2f(0,0)\int_{\mathbb{R}}h(y)
  \exp\{-(\lambda^*+i\omega)(y+c^*\tau)\}dy.
\end{align*}

When $\tau=0$, (\ref{eqV}) reduces to $V_t(\omega,t)=[c_2(\omega)-c_1(\omega)]V(\omega,t)$. Solving this equation yields
\begin{equation*}
  V(\omega,t)=V(\omega,0)\exp\{[c_2(\omega)-c_1(\omega)]t\}.
\end{equation*}
Then we take the inverse Fourier transform,
\begin{equation*}
  \tilde{v}(x,t)=\int_{\mathbb{R}}\tilde{v}(x-\eta,0)
  \mathcal{F}^{-1}[I_0](\eta,t)d\eta,
\end{equation*}
where $I_0(\omega,t)=\exp\{[c_2(\omega)-c_1(\omega)]t\}$. Obviously
\begin{equation}\label{eqV0}
  \begin{split}
    \|\tilde{v}(t)\|_{L^{\infty}(\mathbb{R})}&\leq
    \|\tilde{v}(0)\|_{L^1(\mathbb{R})}\sup\limits_{\eta\in
      \mathbb{R}}\left|\mathcal{F}^{-1}[I_0](\eta,t)d\eta\right|\\
    &\leq
    \frac{1}{2\pi}\|\tilde{v}(0)\|_{L^1(\mathbb{R})}
    \int_{\mathbb{R}}\left|I_0(\omega,t)\right|d\omega.
  \end{split}
\end{equation}
Let $k_1(\omega)=\mathrm{Re}\{c_1(\omega)\}$.It holds that $\left|c_2(\omega)\right|\leq k_2$, where
\begin{equation*}
  k_2= c^*\lambda-d(e^{\lambda^*}+e^{-\lambda^*}-2)-\partial_1f(0,0)
  =\partial_2f(0,0)e^{-\lambda^*c^*\tau}G(\lambda^*)>0.
\end{equation*}
Thus $\left|I_0(\omega,t)\right|\leq \exp\{[k_2-k_1(\omega)]t\}=\exp\{-2td\cos(\lambda^*)(1-\cos\omega)\}$. This, together with (\ref{eqV0}) and Lemma \ref{lem-G}, implies
\begin{equation}\label{eqV01}
  \begin{split}
    \|\tilde{v}(t)\|_{L^{\infty}(\mathbb{R})} &\leq
    \frac{1}{2\pi}\|\tilde{v}(0)\|_{L^1(\mathbb{R})}
    \int_{\mathbb{R}}\exp\{-2td\cos(\lambda^*)(1-\cos\omega)\}d\omega\\
    &\leq C\|\tilde{v}(0)\|_{L^1(\mathbb{R})}t^{-1/2}
  \end{split}
\end{equation}
for all $t>0$. This means the conclusion of Lemma \ref{lem-tiv} holds for $\tau=0$.

For $\tau>0$, we apply Lemma \ref{lem-KH} to equation (\ref{eqV}) and obtain
\begin{equation}\label{eqV2}
  V(\omega,t)= I_1(\omega,t)+\int^0_{-\tau}I_2(\omega,t-s)ds,
\end{equation}
where
\begin{equation*}
  \begin{split}
    I_1(\omega,t)&=
    \exp\{-c_1(\omega)(t+\tau)\}e^{\mathcal{B}(\omega)t}_{\tau}V_0(\omega,-\tau),\\
    I_2(\omega,t-s) &=
    \exp\{-c_1(\omega)(t-s)\}e^{\mathcal{B}(\omega)(t-\tau-s)}_{\tau}
    \left[\frac{d}{ds}V_0(\omega,s)+c_1(\omega)V(\omega,s)\right],\\
    \mathcal{B}(\omega)&=c_2(\omega)\exp\{c_1(\omega)\tau\}.
  \end{split}
\end{equation*}
Then, by taking the inverse Fourier transform of (\ref{eqV2}), we get
\begin{equation}\label{eqV3}
  \tilde{v}(x,t)=
  \mathcal{F}^{-1}[I_1](x,t)+\int^0_{-\tau}\mathcal{F}^{-1}[I_2](x,t-s)ds.
\end{equation}
So
\begin{equation}\label{eqV4}
  \|\tilde{v}(t)\|_{L^{\infty}(\mathbb{R})}\leq
  \frac{1}{2\pi}\|I_1(t)\|_{L^{1}(\mathbb{R})}+\frac{1}{2\pi}
  \int^0_{-\tau}\|I_2(t-s)\|_{L^{1}(\mathbb{R})}ds.
\end{equation}
Note that $\left|\exp\{-c_1(\omega)(t+\tau)\}\right|= \exp\{-k_1(\omega)(t+\tau)\}$. Thus
\begin{equation*}
  \left|\mathcal{B}(\omega)\right| =
  \left|c_2(\omega)\exp\{c_1(\omega)\tau\}\right| \leq
  k_2\exp\{k_1(\omega)\tau\}:=k_3(\omega)
\end{equation*}
and hence
\begin{equation*}
  \left|e^{\mathcal{B}(\omega)t}_{\tau}\right|\leq
  e^{k_3(\omega)t}_{\tau}.
\end{equation*}
Obviously $k_1(\omega)\geq k_1(0)=k_2$. It follows from Lemma \ref{lem-KH} that
\begin{equation}\label{eqE}
  e^{-k_1(\omega)}e^{k_3(\omega)t}_{\tau}\leq
  Ce^{-\varepsilon[k_1(\omega)-k_2]t}=C\exp\{-2td\varepsilon
  \cosh(\lambda^*)(1-\cos\omega)\}
\end{equation}
for some constants $C>0$ and $\varepsilon\in (0,1)$. Therefore, we have
\begin{equation}\label{eqI1}
  \begin{split}
    \|I_1(t)\|_{L^1(\mathbb{R})} &\leq
    C\int_{\mathbb{R}}\exp\{-2td\varepsilon
    \cosh(\lambda^*)(1-\cos\omega)\}
    \left|V_0(\omega,-\tau)\right|d\omega\\
    &\leq  C\int_{\mathbb{R}}\exp\{-2td\varepsilon
    \cosh(\lambda^*)(1-\cos\omega)\} d\omega
    \|\tilde{v}_0(-\tau)\|_{L^1(\mathbb{R})}.
  \end{split}
\end{equation}
This, together with Lemma \ref{lem-G}, implies that
\begin{equation}\label{eqI1-2}
  \|I_1(t)\|_{L^1(\mathbb{R})} \leq
  C\|\tilde{v}_0(-\tau)\|_{L^1(\mathbb{R})}t^{-1/2}.
\end{equation}
Similarly, we have
\begin{equation}\label{eqI2}
  \|I_2(t-s)\|_{L^1(\mathbb{R})}  \leq
  C\left\{\|\partial_t\tilde{v}_0(s)\|_{L^1(\mathbb{R})}
    +\|\tilde{v}_0(s)\|_{L^1(\mathbb{R})}\right\}(t-s)^{-1/2}.
\end{equation}
It follows from (\ref{eqV4}), (\ref{eqI1-2}) and (\ref{eqI2}) that
the conclusion of Lemma \ref{lem-tiv} holds for $\tau>0$. Thus the
proof is complete. \qed

\begin{lemma}
  For $c=c^*$, there holds
  \begin{equation}\label{ch525*}
    \sup\limits_{x\in \mathbb{R}}|u^+(x,t)-\phi(x+ct)|\leq
    Ct^{-1/2},\quad t> 0.
  \end{equation}
\end{lemma}

\proof In view of Lemma \ref{lem-tiv}, we have
\begin{equation*}
  \mbox{$\|\bar{v}(t)\|_{L^{\infty}_{w}(\mathbb{R})}\leq C_1t^{-1/2}$
    for all $t> 0$ and some constant $C_1>0$.}
\end{equation*}
This, together with the fact that $v(\xi,t)\leq \bar{v}(\xi,t)$ for $(\xi,t)\in \mathbb{R}\times \mathbb{R}^+$ and $w(\xi)\geq 1$ for $\xi\in (-\infty,x_0]$, implies that
\begin{equation}\label{import1}
  \mbox{$\|v(t)\|_{L^{\infty}(-\infty,x_0]}\leq C_1t^{-1/2}$ for all
    $t> 0$.}
\end{equation}

In what follows, we show that
\begin{equation}\label{import0}
  \mbox{$\|v(t)\|_{L^{\infty}[x_0,\infty)}\leq Ct^{-1/2}$ for all $t>
    0$ and some constant $C>0$.}
\end{equation}
First we see that $v=v(\xi,t)$ satisfies
\begin{equation}\label{ch51*2-1}
  \begin{cases}
    v_t+c^*v_{\xi}-d\cdot
    \Delta_1v-G_1v-G_2h*v_{\tau}\leq 0 &\mbox{for $t>0$ and $\xi\geq
      x_0$},\\
    v(x_0,t)\leq C_1t^{-1/2} &\mbox{for $t>0$ },\\
    v(\xi,s)\leq v_0(\xi,s)&\mbox{for $s\in [-\tau,0]$ and
      $\xi\in \mathbb{R}$}.
  \end{cases}
\end{equation}
Here $x_0\in \mathbb{R}$ is chosen to be sufficiently large such that $G_1(\xi)+G_2(\xi)<0$ for all $\xi\geq x_0$.

Let
\begin{equation*}
  \bar{v}(\xi,t)=C_2(1+t+\tau)^{-1/2},\quad t\in [t^*,\infty),
\end{equation*}
where $t^*$ and $C_2>v_0(\xi,s)\geq 0$ for all $s\in \mathbb{R}$ are sufficient large constants such that
\begin{multline*}
  \bar{v}_t+c^*\bar{v}_{\xi}-d\cdot
  \Delta_1\bar{v}-G_1\bar{v}-G_2h*\bar{v}_{\tau} \\
  =
  \bar{v}\left[-\frac{1}{2}(1+t+\tau)^{-1}-G_1(\xi)-G_2(\xi)
    \frac{1+t+\tau}{1+t}\right]
  \geq  0
\end{multline*}
for all $t\geq t^*$ and $\xi\geq x_0$. In addition, we can choose $C_2$ large enough such that $\bar{v}(\xi,t)\geq v(\xi,t)$ for all $s\in [0,t_0]$ and $\xi\in \mathbb{R}$. This implies that $\bar{v}(\xi,t)$ is an upper solution of (\ref{ch51*2-1}). Hence
\begin{equation*}
  0\leq v(\xi,t)\leq \bar{v}(\xi,t)\leq Ct^{-1/2}
\end{equation*}
for all $t>0$ and $\xi\geq x_0$. This proves (\ref{import0}).  Combing (\ref{import1}) and (\ref{import0}), we obtain the decay rates for $v(\xi,t)$ in $L^{\infty}(\mathbb{R})$,
\begin{equation*}
  \mbox{$\|v(t)\|_{L^{\infty}(\mathbb{R})}\leq C_1t^{-1/2}$ for all
    $t> 0$.}
\end{equation*}
This, together with the fact $v(\xi,t) = u^+(x,t)-\phi(x+c^*t)$, proves this lemma.  \qed

\step{Step 2. The convergence of $u^-(x,t)$ to $\phi(x+ct)$}

For any $c\geq c^*$, let $\xi=x+ct$ and
\begin{equation*}
  v(\xi,t)=\phi(x+ct)-u^-(x,t),\quad v_0(x,t)=\phi(x+cs)-u^-_0(x,s).
\end{equation*}
It follows that $v(\xi,t)\geq 0$ and $v_0(\xi,s)\geq 0$. By using similar arguments as before, we can prove that $u^-(x,t)$ converges to $\phi(x + ct)$ as follows.

\begin{lemma}
  There holds the exponential decay
  \begin{equation}\label{ch533}
    \sup\limits_{x\in \mathbb{R}}|\phi(x+ct)-u^-(x,t)|\leq
    Ce^{-\mu t},\quad t\geq 0,\quad \mbox{for $c>c^*$},
  \end{equation}
  and the algebraic decay
  \begin{equation}\label{ch533*}
    \sup\limits_{x\in \mathbb{R}}|\phi(x+ct)-u^-(x,t)|\leq
    C(1+t)^{-1/2},\quad t\geq 0,\quad \mbox{for $c=c^*$},
  \end{equation}
  where $\mu>0$ satisfies $\mathcal{M}(c,3\mu)>0$ and   $\mathcal{N}(3\mu)<0$.
\end{lemma}

\step{Step 3. The convergence of $u(x,t)$ to $\phi(x+ct)$}

Finally we prove that $u(x,t)$ converges to $\phi(x+ct)$. In fact, notice that $u^-(x,t)\leq u(x,t)\leq u^+(x,t)$ for all $(x,t)\in \mathbb{R}\times [0,\infty)$ and combine (\ref{ch525}), (\ref{ch525*}), (\ref{ch533}), and (\ref{ch533*}). This yields the following theorem.

\begin{theorem}
  There hold the exponential decay
  \begin{equation*}
    \sup\limits_{x\in \mathbb{R}}|u(x,t)-\phi(x+ct)|\leq
    Ce^{-\mu t},\quad t\geq 0,\quad \mbox{for $c>c^*$},
  \end{equation*}
  and the algebraic decay
  \begin{equation*}
    \sup\limits_{x\in \mathbb{R}}|u(x,t)-\phi(x+ct)|\leq
    C(1+t)^{-1/2},\quad t\geq 0,\quad \mbox{for $c=c^*$},
  \end{equation*}
  where $\mu>0$ satisfies $\mathcal{M}(c,3\mu)>0$ and   $\mathcal{N}(3\mu)<0$.
\end{theorem}

\section{Applications}
\label{sec:Applications}

In this section, we give applications of the stability result in combination with the results of \cite{Guo:11a} to some biological and epidemiological models.

\subsection{A host-vector disease model}

We consider the host-vector disease model with discrete delay
\begin{equation}\label{ex1}
  u_t(x,t)=d\cdot\Delta_1u(x,t)-au(x,t)+bu(x,t-\tau)[1-u(x,t)],
  \quad   x\in \mathbb{R},\,t>0,
\end{equation}
where the function $u(x, t)$ denotes the density of the infectious host at time $t > 0$ and spatial location $x\in \mathbb{R}$, $\tau\geq 0$ stands for the time delay; $b>a\geq 0$ are parameters.  For the continuous host-vector disease model, i.e., (\ref{ex1}) with $\Delta_1u(x,t)$ replaced by $u_{xx}(x,t)$, the existence of traveling wave was shown by Lin and Hong \cite{LinHong}, and later the nonlinear stability of noncritical waves was investigated by Lv and Wang \cite{LvWang}.

Obviously, (\ref{ex1}) has two equilibria $u= 0$ and $K=1-\frac{a}{b}$. We can rewrite (\ref{ex1}) in the form of (\ref{eq}) with $f(u,v)=-au+bv(1-u)$ and $h(x)=\delta(x)$, where $\delta(\cdot)$ is the Dirac delta function. Obviously, $f(0,0)=f(K,K)=0$ and $f(u,u)>0$ for all $u\in (0,K)$. Moreover, it is straightforward to verify that the
conditions (F1)--(F3) hold. Thus we have the following result.

\begin{theorem}\label{thm-vd} There exists a
  minimal wave speed $c^*>0$ such that for each $c\geq c^*$,   (\ref{ex1}) has exactly one traveling front $\phi(x+ct)$ (up to  translation) connecting $u=0$ and $u=1-\frac{a}{b}$. Moreover, the   wave profile $\phi$ has the following properties.
  \begin{description}
  \item[(i)] $0<\phi(\xi)<1-\frac{a}{b}$ and $\phi'(\xi)>0$ for all     $\xi\in \mathbb{R}$.
  \item[(ii)] Both $\lim_{\xi\to-\infty}\phi'(\xi)/\phi(\xi)$ and     $\lim_{\xi\to\infty}\phi'(\xi)/[b-a-b\phi(\xi)]$ exist and are     positive.
  \item[(iii)] The wave front $\phi(x + ct)$ is nonlinearly stable in     the sense of Theorem \ref{thm-sta}. Moreover, the exponential     decay rate $\mu>0$ is less than the minimal positive root of     $3\mu+ae^{3\mu\tau}-b=0$.
  \end{description}
\end{theorem}

\begin{remark}
  Taking $\tau=0$, $a=0$, and $b=1$ in (\ref{ex1}), we obtain the   following well-known Fisher-KPP equation
  \begin{equation}\label{ex4}
    u_t(x,t)=d\cdot\Delta_1u(x,t)-u(x,t)[1-u(x,t)],\quad x\in
    \mathbb{R},\,t>0.
  \end{equation}
  Thus, Theorem \ref{thm-vd} can be applied to (\ref{ex4}).
\end{remark}

\subsection{An nonlocal population model with age structure}

Let $h(x)=\frac{1}{\sqrt{4\pi\alpha}}\exp\{-\frac{y^2}{4\alpha}\}$ and $f(u,v)=-\delta u^2+pe^{-\gamma\tau}v$ with $\alpha>0$, $\delta>0$, $p>0$ and $\gamma>0$. We then reduce (\ref{eq}) to the age-structured population model
\begin{equation}\label{ex2}
  u_t(x,t)=d\cdot\Delta_1u(x,t)-\delta
  u^2(x,t)+pe^{-\gamma\tau}\int_{\mathbb{R}}g_{\alpha}(y)u(x-y,t-\tau)dy,\quad
  x\in \mathbb{R},\,t>0.
\end{equation}
It is clear that the constant equilibria of (\ref{ex2}) are $u=0$ and $u=\frac{p}{\delta}e^{-\gamma\tau}$. Moreover, $\partial_1f(u,v)=-2\delta u\leq 0$ and $\partial_2f(u,v)=pe^{-\gamma\tau}$. Obviously, conditions (F1)--(F3) and (H1)--(H2) are satisfied. Then the result of \cite{Guo:11a} shows that there exists a minimal wave speed $c^*>0$ such that for each $c\geq c^*$, (\ref{ex2}) has exactly one (up to translation) monotone traveling front $\phi(x+ct)$ (up to a translation) connecting $u=0$ and $u=\frac{p}{\delta}e^{-\gamma\tau}$. The stability result the is as follows.
\begin{proposition}
  The wave front $\phi(x + ct)$ is nonlinearly stable in the sense of   Theorem \ref{thm-sta}. Moreover, the exponential decay rate $\mu>0$   is less than the minimal positive root of   $3\mu+pe^{-\gamma\tau}[e^{3\mu\tau}-2]=0$.
\end{proposition}

\begin{remark}
  For the spatially continuous age-structured population model, i.e.,   (\ref{ex2}) with $\Delta_1u(x,t)$ replaced by $u_{xx}(x,t)$, which   was first derived in \cite{Al-Omari2}, the linear and nonlinear   stability of noncritical traveling waves are shown in   \cite{Gourley2005} and \cite{LiMei}, respectively, in the case where   the time delay $\tau$ is sufficiently small. This restriction was   removed in \cite{MeiWong}.
\end{remark}

\subsection{A nonlocal Nicholson's blowflies model with delay}

Let $f(u,v)=-\delta u+pve^{-av}$ with $a>0$ and $p>\delta>0$. We then reduce (\ref{eq}) to the delay model
\begin{multline}\label{ex3}
  u_t(x,t)=d\cdot\Delta_1u(x,t)-\delta
  u(x,t)\\
  +p(h*u)(x,t-\tau)\exp\{-a(h*u)(x,t-\tau)\},
  \quad x\in \mathbb{R},\,t>0.
\end{multline}
It is easy to see that (\ref{ex3}) has two equilibria $u=0$ and $u=\frac{1}{a}\ln\frac{p}{\delta}$. Equation (\ref{ex3}) is a discrete analog of the following nonlocal diffusive Nicholson's blowflies equation with delay (see \cite{Gourley-Ruan}), for $x\in   \mathbb{R}$ and $t>0$,
\begin{equation}\label{ex3'}
  u_t(x,t)=du_{xx}(x,t)-\delta
  u(x,t)+p(h*u)(x,t-\tau)\exp\{-a(h*u)(x,t-\tau)\},
\end{equation}
Recently, some special cases of (\ref{ex3'}) has been investigated in \cite{Gourley,Gurney,Lin, Mei1, Mei2, So2}.  In particular, Mei \emph{et al}. studied the stability of traveling wave fronts of (\ref{ex3'}) with $h$ being the Dirac delta function, namely,
\begin{equation}\label{ex3''}
  u_t(x,t)=du_{xx}(x,t)-\delta
  u(x,t)+pu(x,t-\tau)\exp\{-au(x,t-\tau)\},
  \quad x\in \mathbb{R},\,t>0.
\end{equation}

Similar to the two previous examples, there exists for a minimal wave speed $c^*>0$ such that there is a unique (up to translation) monotone traveling wave solution to (\ref{ex3}) connecting $u=0$ and $u=\frac{1}{a}\ln\frac{p}{\delta}$, provided that $c\geq c^*$ \cite{Guo:11a}.
\begin{proposition}
  Assume further that $\delta(\ln p-\ln\delta)<2\delta-p$, then the   wave front $\phi(x + ct)$ is nonlinearly stable in the sense of   Theorem \ref{thm-sta}. Moreover, the exponential decay rate $\mu>0$   is less than the minimal positive root of $3\mu-\delta-\delta   e^{3\mu\tau}(1+\ln \delta-\ln p)=0$.
\end{proposition}

\begin {thebibliography}{99}

\bibitem{Al-Omari2} J. Al-Omari and S. A. Gourley, A nonlocal reaction-diffusion
model for a single species with stage structure and distributed
maturation delay, European J. Appl. Math. 70 (2005) 858--879.

\bibitem{Al} J.~Al-Omari and S.A.~Gourley, Monotone traveling fronts in age-structured reaction-diffusion model of a
single species, J. Math. Biol. 45  (2002) 294--312.

\bibitem{Bates1} P.W. Bates, P.C. Fife, X.F. Ren and X.F. Wang, Traveling waves in a convolution model for phase
transitions, Arch. Rational Mech. Anal. 138 (1997) 105--136. 

\bibitem{Bates2} P.W. Bates and A. Chmaj, A discrete convolution model for phase transitions, Arch. Rational Mech.
Anal. 150 (1999) 281--305. 

\bibitem{Bell} J. Bell, Some threshold results for models of myelinated nerves, Math. Biosci. 54 (1981) 181--190.


\bibitem{Berestycky} H. Berestycki, G.  Nadin, B. Perthame, L. Ryzhik, The non-local Fisher-KPP equation: travelling waves and steady states. Nonlinearity 22 (2009), no. 12, 2813--2844.

\bibitem{Britton1} N. Britton, Aggregation and the competitive exclusion principle, J. Theor. Biol. 136 (1989), 57--66.

\bibitem{Britton2} N. Britton, Spatial structures and periodic travelling waves in an integro-differential reaction-diffusion population model, SIAM J. Appl. Math. 50 (1990), 1663--1688.

\bibitem{Cahn} J.W. Cahn, J. Mallet-Paret and E.S. van Vleck, Traveling wave solutions for systems of ODEs on a two-dimensional spatial lattice, SIAM J. Appl. Math. 59 (1999) 455--493. 

\bibitem{Carpio} A. Carpio and G. Duro, Explosive behavior in spatially discrete reaction-diffusion systems, Discrete Contin. Dyn. Syst. Ser. B 12 (2009), 693--711.

\bibitem{Chen1997} X. Chen,  Existence, uniqueness, and asymptotic stability
of traveling waves in nonlocal evolution equations,  Adv.
Differential Equations 2 (1997) 125--160.

\bibitem{Chua} L.O. Chua, T. Roska, The CNN paradigm, IEEE Trans. Circuits Syst. 40 (1993) 147--156.

\bibitem{Coville} J. Coville, L. Dupaigne, On a nonlocal reaction diffusion equation arising in population dynamics,
Proceedings of the Royal Society of Edinburgh 137  (2007) A, 1--29.

\bibitem{Fisher} R.A. Fisher, The advance of advantageous genes, Ann. Eugenics 7 (1937) 355--369.

\bibitem{Gallay} T. Gallay, Local stability of critical fronts in nonlinear parabolic PDEs. Nonlinearity 7
(1994) 741--764.

\bibitem{Chow} S.-N. Chow, J. Mallet-Paret and W. Shen, Traveling waves in lattice dynamical systems, J. Differential Equations 149 (1998) 248--291. 

\bibitem{Fath} G. F\'ath, Propagation failure of traveling waves in a discrete bistable medium, Phys. D 116 (1998) 176--190.

\bibitem{Gourley2005} S. Gourley, Linear stability of traveling fronts in an age-structured reaction-diffusion population model. Q. J. Mech.
Appl. Math. 58 (2005) 257--268.

\bibitem{Gourley0} S.A. Gourley, J.W.-H. So and J.H. Wu, Non-locality of reaction-diffusion equations induced by delay:
biological modelling and nonlinear dynamics, in: D.V. Anosov, A. Skubachevskii (Guest Eds.),
Contemporary Mathematics. Thematic Surveys, Kluwer Plenum, Dordrecht, New York, 2003, pp.
84--120.


\bibitem{Gourley} S. A. Gourley and J. Wu, Delayed non-local diffusive systems in biological invasion and
disease spread, in Nonlinear Dynamics and Evolution Equations, Fields Inst. Commun. 48,
H. Brunner, X.-Q. Zhao, and X. Zou, eds., AMS, Providence, RI, 2006, pp. 137--200.


\bibitem{Gourley-Ruan} S.A. Gourley, S. Ruan, Dynamics of the diffusive Nicholson's blowflies equation with distributed delay, Proc. Roy. Soc. Edinburgh Sect. A 130 (2000) 1275--1291.

\bibitem{Guo:11a} S. Guo and J. Zimmer, Travelling wavefronts in discrete reaction-diffusion equations   with nonlocal delay effects, arXiv:1406.5321 [math.DS].

\bibitem{Gurney} W. S. C. Gurney, S. P. Blythe, and R. M. Nisbet, Nicholson's blowflies revisited, Nature, 287 (1980) 17--21.

\bibitem{Hale} J.K. Hale, S.M. Verduyn Lunel, Introduction to functional-differential equations. Applied Mathematical Sciences 99. Springer-Verlag, New York, 1993


\bibitem{Harterich:02a} J. H{\"a}rterich, B. Sandstede, A. Scheel, Exponential dichotomies for linear non-autonomous functional differential equations of mixed type, Indiana Univ. Math. J. 51 (2002) 1081--1109.

\bibitem{Huang} R. Huang, M. Mei, and Y. Wang, Planar traveling waves for nonlocal dispersal equation with monostable nonlinearity,  Discrete Contin. Dyn. Syst. -- Series A, Vol. 32 (2012) 3621--3649.

\bibitem{Hupkes09a} H. J. Hupkes, B. Sandstede, Modulated wave trains in lattice differential systems, J. Dynam. Differential Equations 21 (2009) 417--485.

\bibitem{Keener} J.P. Keener, Propagation and its failure in coupled systems of discrete excitable cells, SIAM J. Appl.
Math. 22 (1987) 556--572.

\bibitem{Kirchgassner} K. Kirchg\"assner, On the nonlinear dynamics of traveling fronts. J. Diff. Eqns 96 (1992)
256--278.

\bibitem{Kolmogorov} A.N. Kolmogorov, I.G. Petrovsky and N.S. Piskunov, \'Etude de l'\'equation de la diffusion avec
croissance de la quantit\'e de mati\'ere et son application \'a un probl\'eme biologique, Bull. Univ.
Moskov. Ser. Internat., Sect. A 1 (1937) 1--25

\bibitem{Laplante-Erneux} J.P. Laplante, T. Erneux, Propagation failure in arrays of coupled bistable chemical reactors, J. Phys. Chem. 96 (1992) 4931--4934.

\bibitem{LiMei} G. Li, M. Mei, Y. Wong, Nonlinear stability of traveling wavefronts in an age-structured reaction-diffusion
population model. Math. Biosci. Eng. 5  (2008) 85--100.





\bibitem{Liang1} X. Liang and X.-Q. Zhao, Asymptotic speeds of spread and traveling waves for monotone
semiflows with applications, Comm. Pure Appl. Math. 60 (2007) 1--40.


\bibitem{Liang2} X. Liang and X.-Q. Zhao, Erratum: Asymptotic speeds of spread and traveling waves for
monotone semiflows with applications, Comm. Pure Appl. Math. 61 (2008) 137--138.



\bibitem{LinHong} G. Lin and Y. Hong, Travelling wavefronts in a vector disease
model with delay, Appl. Math. Modelling 32 (2008) 2831--2838.

\bibitem{Lin} C.-K. Lin and M. Mei, On travelling wavefronts of the Nicholson's blowflies equations with diffusion, Proc. Roy. Soc. Edinburgh Sect. A, 140 (2010) 135--152.

\bibitem{LvWang} G. Lv and M. Wang, Nonlinear stability of travelling wave fronts
for delayed reaction diffusion equations, Nonlinearity 23 (2010)
84--873.


\bibitem{LvWang2010} G. Lv and M. Wang, Existence, uniqueness and stability of
traveling wave fronts of discrete quasi-linear equation, Discrete and Continuous Dynamical System B 13 (2010) 415--433.

\bibitem{LvWang2012} G. Lv and M. Wang, Nonlinear stability of traveling wave fronts for nonlocal delayed reaction-diffusion
equations. J. Math. Anal. Appl. 385 (2012) 1094--1106.



 \bibitem{Lucia} M. Lucia,
C. Muratov, and  M. Novaga, Linear vs. nonlinear selection for the
propagation speed of the solutions of scalar reaction-diffusion equations
invading an unstable equilibrium, Comm. Pure Appl. Math. 57 (2004)
616--636.

\bibitem{Ma} S. Ma and X. Zou, Existence, uniqueness and stability of travelling
waves in a discrete reaction-diffusion monostable
equation with delay, J. Differential Equations 217 (2005) 54--87.



\bibitem{Mei1} M. Mei, C.-K. Lin, C.-T. Lin, and J. W.-H. So, Traveling wavefronts for time-delayed
reaction-diffusion equation: (I) local nonlinearity, J. Differential Equations, 247 (2009) 495--510.

\bibitem{Mei4} M. Mei, C. Ou, X. Zhao, Global stability of monostable traveling waves for
nonlocal time-delayed reaction-diffusion equations, SIAM J. Math. Anal.
42 (2010) 2762--2790.

\bibitem{Mei4a} M. Mei, C. Ou, X. Zhao, Erratum: Global stability of monostable traveling waves for nonlocal time-delayed reaction-diffusion equations. SIAM J. Math. Anal.  44  (2012), 538--540

\bibitem{Mei3} M. Mei, J.W.-H. So, Stability of strong traveling waves for a nonlocal time-delayed reaction-diffusion equation. Proc.
Roy. Soc. Edinb. 138A (2008) 551--568.


\bibitem{Mei2} M. Mei, J.W.-H. So, M. Li, and S. Shen, Asymptotic stability of t
raveling waves for the Nicholson's blowflies equation with diffusion, Proc. Roy. Soc. Edinburgh Sect. A, 134 (2004)
579--594.

\bibitem{Mei-Wang} M. Mei and Y. Wang, Remark on stability of planar traveling waves
for nonlocal Fisher-KPP type reaction-diffusion equations in
multi-dimensional space, International Journal Of Numerical Analysis
And Modeling--Series B, in press.


\bibitem{MeiWong} M. Mei and Y. S. Wong, Novel stability results for traveling
wavefronts in an age-structured reaction-diffusion equations, Math.
Biosci. Eng. 6 (2009) 743--752.

\bibitem{Ou} C. Ou, J. Wu, Persistence of wavefronts in delayed nonlocal reaction-diffusion equations. J. Differential Equations  235  (2007) 219--261

\bibitem{Ruan} S. Ruan and D. Xiao, Stability of steady states and existence of traveling waves in a vector disease
model. Proc. Roy. Soc. Edinburgh 134 (2004) 991--1011 .

\bibitem{Rustichini1} A.  Rustichini,  Functional  differential  equations  of  mixed  type:  The  linear  autonomous  case,
J.  Dynam.  Differential  Equations  1  (1989)  121-143.

\bibitem{Rustichini2} A.  Rustichini,  Hopf  bifurcation  for  functional  differential  equations  of  mixed  type,  J.  Dynam.
Differential  Equations  1  (1989) 145-177.

\bibitem{Sattinger} D. H. Sattinger, On the stability of waves of nonlinear parabolic systems. Adv. Math. 22
(1976), 312--355.

\bibitem{sco} A. Scott, Nonliear Science: Emergence and dynamics of coherent structures. Oxford Texts in Applied and Engineering Mathematics, 8. Oxford University Press 2003.

\bibitem{Smith} H.L. Smith and H. Thieme, Strongly order preserving semiflows generated by functional differential
equations, J. Differential Equations 93 (1991) 332--363.

\bibitem{SmithZhao} H. L. Smith and X.-Q. Zhao, Global asymptotic stability of travelling waves in delayed
reaction-diffusion equations. SIAM J. Math. Analysis 31 (2000) 514--534.

\bibitem{Khusainov} D. Ya Khusainov, A. F. Ivanov and I. V. Kovarzh, Solution of one
heat equation with delay, Nonlinear Oscillations, 12 (2009)
260--282.







\bibitem{Schaaf} K. W. Schaaf, Asymptotic behavior and traveling wave solutions for parabolic functional
differential equations. Trans. Am. Math. Soc. 302 (1987) 587--615.

\bibitem{So} J. W.-H. So, J. Wu and X. Zou, A reaction-diffusion model for a single species with age
structure: (I) Traveling wavefronts on unbounded domains, Proc. R. Soc. Lond. Ser. A
Math. Phys. Eng. Sci., 457 (2001) 1841--1853. 

\bibitem{So2} J. W.-H. So and X. Zou, Traveling waves for the diffusive Nicholson's blowflies equation,
Appl. Math. Comput., 22 (2001) 385--392.

\bibitem{Thieme} H. Thieme and X.-Q. Zhao, Asymptotic speeds of spread and traveling waves for integral equation
and delayed reaction-diffusion models, J. Differential Equations, 195 (2003) 430--370.

\bibitem{Wang} Z. Wang, W. Li, S. Ruan, Travelling fronts in monostable equations with nonlocal delayed effects. J. Dyn. Differ.
Equ. 20 (2008) 573--607

\bibitem{WuZou2} J. Wu, X. Zou, Local existence and stability of periodic traveling waves of lattice functional differential
equations, Canad. Appl. Math. Quart. 6 (1998) 397--416.






\end {thebibliography}
\end{document}